\documentclass {article}
\usepackage{authblk}
\usepackage{etex}
\usepackage[utf8]{inputenc}
\usepackage[english]{babel}
\usepackage{amsmath}
\usepackage{amsthm}
\usepackage{amssymb}
\usepackage{sectsty}
\usepackage{titlesec}
\usepackage{color}
\usepackage[color,matrix,arrow]{xy}
\usepackage{amsgen}
\usepackage{amstext}
\usepackage{amsbsy}
\usepackage{amsopn}
\usepackage{amsfonts}
\usepackage{eepic}
\usepackage{graphicx}
\usepackage{epsf}
\usepackage{pstricks}
\usepackage{enumerate}
\xyoption{all}

\usepackage{pgf}
\usepackage{tikz}

\newcommand{\footrecall}[1]{}

\titleformat*{\section}{\large\bfseries}
\titleformat*{\subsection}{\normalsize \bfseries}

\title{On the dynamics of extensions of free-abelian times free groups endomorphisms to the completion}

\author{Andr\'e Carvalho \thanks{andrecruzcarvalho@gmail.com}}
\affil{Centre of Mathematics, University of Porto, R. Campo Alegre, 4169-007 Porto,
Portugal}

\date{}

\begin{document}

\newtheorem{theorem}{Theorem}[section]
\newtheorem{lemma}[theorem]{Lemma}
\newtheorem{question}[theorem]{Question}
\newtheorem{remark}[theorem]{Remark}
\newtheorem*{theorem*}{Theorem}
\newtheorem{proposition}[theorem]{Proposition}
\newtheorem{corollary}[theorem]{Corollary}

\theoremstyle{definition}
\newtheorem{definition}[theorem]{Definition}
\newtheorem{example}[theorem]{Example}

\newcommand{\ophi}{\overline{\varphi}}
\newcommand{\opsi}{\overline{\psi}}
\newcommand{\N}{\mathbb{N}}
\newcommand{\Z}{\mathbb{Z}}
\newcommand{\F}{\mathbb{F}}
\newcommand{\R}{\mathbb{R}}
\newcommand{\Q}{\mathbb{Q}}
\newcommand{\Fix}{\text{Fix}}
\newcommand{\Sing}{\text{Sing}}
\newcommand{\Per}{\text{Per}}
\newcommand{\Ker}{\text{Ker}}
\newcommand{\End}{\text{End}}
\newcommand{\Mono}{\text{Mon}}
\newcommand{\Aut}{\text{Aut}}
\newcommand{\Reg}{\text{Reg}}
\newcommand{\mc}{\mathcal}

\maketitle

\begin{abstract}
We obtain conditions of uniform continuity for endomorphisms of free-abelian times free groups for the product metric defined by taking the prefix metric in each component and establish an equivalence between uniform continuity for this metric and the preservation of a coarse-median, which was recently introduced by Fioravanti. Considering the extension of an endomorphism to the completion we count the number of orbits for the action of the subgroup of fixed points (resp. periodic) points on the set of infinite fixed (resp. periodic) points. Finally, we study the dynamics of infinite points: for some  endomorphisms, defined in a precise way, fitting a classification given by Delgado and Ventura, we prove that every infinite point is either periodic or wandering, which implies that the dynamics is asymptotically periodic. We also prove the latter for the case of automorphisms.
\end{abstract}

\section{Introduction}

The dynamical study of endomorphisms of groups started with the (independent) work of Gersten \cite{[Ger87]} and Cooper \cite{[Coo87]}, using respectively graph-theoretic and topological approaches. They proved that the subgroup of fixed points $\Fix(\varphi)$ of some fixed automorphism $\varphi$ of $F_n$ is always finitely generated, and Cooper succeeded on classifying from the dynamical viewpoint the fixed points of the continuous extension of $\varphi$ to the boundary of $F_n$. Bestvina and Handel subsequently developed the theory of train tracks to prove that $\Fix(\varphi)$ has rank at most $n$ in \cite{[BH92]}. The problem of computing a basis for $\Fix(\varphi)$ had a tribulated history and was finally settled by Bogopolski and Maslakova in 2016 in \cite{[BM16]}. 

This line of research extended early to wider classes of groups. For instance, Paulin proved in 1989 that the subgroup of fixed points of an automorphism of a hyperbolic group is finitely generated \cite{[Pau89]}. Fixed points were also studied for right-angled Artin groups \cite{[RSS13]} and lamplighter groups \cite{[MS18]}. 

Regarding the extension of an endomorphism to the completion, infinite fixed points of automorphisms of free groups were also discussed by Bestvina and Handel in  \cite{[BH92]} and Gaboriau, Jaeger, Levitt and Lustig in \cite{[GJLL98]}. The dynamics of free groups automorphisms is proved to be asymptotically periodic in \cite{[LL08]}. In \cite{[CS09a]}, Cassaigne and Silva study the dynamics of infinite fixed points for monoids defined by special confluent rewriting systems (which contain free groups as a particular case). This was also achieved by Silva in \cite{[Sil13]} for virtually injective endomorphisms of virtually free groups.

In this paper we will study the dynamics of infinite points of free-abelian times free groups which is a subclass of the well-known right-angled Artin groups. This class of groups have been thoroughly studied both algebraically and algorithmically in \cite{[DV13]}.
There, Delgado and Ventura prove that for $G=\mathbb Z^m\times F_n$, with $n\neq 1$, all endomorphisms of $G$ are of one of the following forms: 
\begin{enumerate}[(I)]
\item $\Psi_{\Phi,{Q, P}}:(a,u) \mapsto(aQ+\mathbf uP, u\Phi)$, where $\Phi\in End(F_n)$, $ Q\in \mc M_m(\mathbb Z)$, and $ P\in \mc M_{n\times m}(\mathbb Z)$.
\item $\Psi_{z,{\ell,h,Q,P}}:(a,u)\mapsto (aQ+\mathbf uP,z^{{a\ell}^T+{\mathbf uh}^T})$, where $1\neq z\in F_n$ is not a proper power, $ Q\in \mc M_m(\mathbb Z)$, $ P\in \mc M_{n\times m}(\mathbb Z)$, $\mathbf 0\neq \ell\in \mathbb Z^m$, and $ h\in\mathbb Z^n$,
\end{enumerate}
where $\mathbf u\in \mathbb Z^n$ denotes the abelianization of the word $u\in F_n$. This classification will be a very important tool for us.
 Although the class of free-abelian groups is very well known and the class of free groups, while being much more complex than the first, has also been deeply studied, some problems in the product $\mathbb Z^m\times F_n$ are not easily reduced to problems in each factor. In particular, when endomorphisms (and automorphisms) are considered, we have that many endomorphisms are not obtained by applying an endomorphism of $\mathbb Z^m$ to the first component and one of $F_n$ to the second, so some problems arise when the dynamics of an endomorphism is considered.

We denote by $F_n$ the free group of rank $n$ and its alphabet by $A=\{a_1,\ldots,a_n\}$. Given two words $u$ and $v$ on a free group, we write $u\wedge v$ to denote the longest common prefix of $u$ and $v$. The prefix metric on a free group is defined by
$$d(u,v)=\begin{cases}
2^{-|u\wedge v|} \text{ if $u\neq v$}\\
0 \text{ otherwise}
\end{cases}.$$
The prefix metric on a free group is in fact an ultrametric and its completion $(\hat{F}_n,\hat{d})$ is a compact space which can be described as the set of all finite and infinite reduced words on the alphabet $A \cup A^{-1}$. We will denote by $\partial F_n$ the set consisting of only the infinite words and call it the \emph{boundary} of $F_n$.

A free-abelian times free group is of the form $\mathbb Z^m\times F_n$, which we consider endowed with the product metric given by taking the prefix metric in each (free) component, i.e.,
 $$d((a,u),(b,v))=\max\{d(a_1,b_1),\ldots, d(a_m,b_m), d(u,v)\},$$
where $a_i$ and $b_i$ denote the $i$-th component of $a$ and $b$, respectively.  This metric is also an ultrametric and $\widehat {\mathbb Z^m \times F_n}$ is homeomorphic to $\hat {\mathbb Z^m}\times \hat {F_n}$ by uniqueness of the completion (Theorem 24.4 in \cite{[SW70]}).

There are ways to think of infinity for this class of groups like the CAT(0) boundary or the Roller boundary. However, this feels like a very natural way to define infinity for this class and the dynamics at the infinity was not previously considered for these groups with this particular definition of boundary. 

It is well known by a general topology result \cite[Section XIV.6]{[Dug78]} that every uniformly continuous mapping $\varphi$ between metric spaces admits a unique continuous extension $\hat\varphi$ to the completion. The converse is obviously true in this case: if a mapping between metric spaces admits a continuous extension to the completion, since the completion is compact, then the extension must be uniformly continuous, and so does the restriction to the original mapping. For this reason, we will focus our work on the uniformly continuous endomorphisms and their extensions to the completion. It is well known that for a free group $F_n$ endowed with the prefix metric, an endomorphism $\varphi\in End(F_n)$ is uniformly continuous if and only if it is either constant or injective. Characterizing and studying some properties of uniformly continuous endomorphisms has been done before for other classes of groups (see for example \cite{[CS09b]}, \cite{[Sil10]}, \cite{[Sil13]}, \cite{[AS16]}, \cite{[Car21]}).

An interesting property of this metric (and so, of this boundary) is that the uniformly continuous endomorphisms of $\mathbb Z^m\times F_n$ for this metric $d$ are precisely the coarse-median preserving endomorphisms for the product coarse median obtained by taking the median operator induced by the metric $\ell_1$ in $\mathbb Z^m$ and the one given by hyperbolicity of $F_n$:

\newtheorem*{cmpfatf}{Theorem \ref{cmpfatf}}
\begin{cmpfatf}
An endomorphism $\varphi\in \End(\mathbb Z^m\times F_n)$ is uniformly continuous with respect to the product metric $d$ obtained by taking the prefix metric in each direct factor if and only if it is coarse-median preserving for the product coarse median $\mu$ obtained by taking the median operator $\mu_1$ induced by the metric $\ell_1$ in $\mathbb Z^m$ and the coarse median operator $\mu_2$ given by hyperbolicity in $F_n$.
\end{cmpfatf}

Coarse-median preservation turns out to be a useful tool to obtain interesting properties of automorphisms (see \cite{[Fio21]}), including finiteness results on the fixed subgroup of an automorphism.

In \cite{[Sil13]}, the author, considering automorphisms of virtually free groups, proves that infinite fixed points that belong to the topological closure $(\Fix(\varphi))^c$ are never attractors nor repellers, while infinite fixed points not belonging to  $(\Fix(\varphi))^c$ must always be either attractors or repellers, generalizing a previously known result for free groups (see \cite{[Coo87]}). 

Fixed points act naturally on infinite fixed points by left multiplication. It is known that, for free groups, the number of $(\Fix(\varphi))$-orbits of the set of attracting fixed points is finite (see \cite{[Coo87]}) and that is used to define the index of an automorphism in \cite{[GJLL98]}. A similar result is obtained for virtually free groups in \cite{[Sil13]}. As highlighted in Section \ref{finitinffixed}, a similar result cannot be obtained in general in this case. However, when $\varphi$ is a type II uniformly continuous endomorphism, then infinite fixed (resp. periodic) points have  $\sum\limits_{i=0}^m 2^i {m\choose i}$ $\Fix(\varphi)$-orbits (resp.  $\Per(\varphi)$-orbits). This result can be seen as some sort of infinite version of Proposition 6.2 in \cite{[DV13]}.

We will also classify infinite fixed points of a uniformly continuous automorphism of $\mathbb Z^m\times F_n$ in attractors or repellers, using the classification obtained for free groups:

\newtheorem*{classinffix}{Corollary \ref{classinffix}}
\begin{classinffix}
Let $\varphi\in Aut(\mathbb Z^m\times F_n)$ be a uniformly continuous automorphism such that $(a,u)\hat\varphi=(a\hat\varphi_1,b\hat\phi)$, where $\varphi_1$ is given by a uniform matrix and $\phi\in Aut(F_n)$. Then an infinite fixed point $(a,u)\in \Fix(\hat\varphi)\setminus \Fix(\varphi)$ is an attractor (resp. repeller) if and only if $a\in Fix(\varphi_1)$ and $u$ is an attractor (resp. repeller) for $\hat\phi$.
\end{classinffix}

In the study of dynamical systems, the notion of \emph{$\omega$-limit set} plays a crucial role. Given a metric space $X$, a continuous function $f:X\to X$, and a point $x\in X$, the $\omega$-limit set $\omega(x,f)$ of $x$ consists of the accumulation points of the sequence of points in the orbit of $x$. Understanding the $\omega$-limits gives us a grasp on the behaviour of the system in the long term. If the space $X$ is compact, then $\omega$-limit sets are nonempty, compact and $f$-invariant.

In \cite{[LL08]}, the authors proved that in the case where $f$ is the extension of a free group automorphism to the completion, then, for every point $x\in \hat F_n$, $\omega(x,f)$ is a periodic orbit.

We will prove that uniformly continuous automorphisms (i.e., the ones that extend to the completion) of $\mathbb Z^m \times F_n$ have the same property:

\newtheorem*{asympdyn}{Theorem \ref{asympdyn}}
\begin{asympdyn}
Every uniformly continuous automorphism $\varphi\in Aut(\mathbb Z^m \times F_n)$ has asymptotically periodic dynamics on $\widehat{\mathbb Z^m\times F_n}$. 
\end{asympdyn}

In the case of uniformly continuous type II endomorphisms, we prove something stronger than this. Informally, a point is said to be \emph{wandering} if it has some neighborhood such that, from some point on, its points leave the neighborhood and don't come back. Obviously, a wandering point cannot belong to an $\omega$-limit set. We will prove that, for a type II uniformly continuous endomorphism, every point in the completion must be either periodic or wandering, showing that non-periodic points, when iterated long enough wander away from some neighborhood \emph{carrying} the neighborhood with them. In particular, since $\omega$-limit sets are nonempty, they must be periodic orbits. In Section \ref{dyninfpoints}, we prove the following theorem:

\newtheorem*{typeiidych}{Theorem \ref{typeiidych}}
\begin{typeiidych}
Let $\varphi\in End(\mathbb Z^m \times F_n)$ be a type II uniformly continuous endomorphism  and consider $\hat\varphi$, its continuous extension to the completion. 
Then every point $(a,u)\in \widehat{\mathbb Z^m \times F_n}$ is either wandering or periodic.
\end{typeiidych}

We also prove that this dichotomy holds for eventually length-nondecreasing automorphisms and some weaker version of it for general eventually length-nondecreasing type I endomorphisms.

The paper is organized as follows.
We find equivalent conditions to uniform continuity of endomorphisms in Section \ref{unicontsec}.  In Section \ref{cmpendos} we establish an equivalence between uniform continuity and preservation a coarse median for  endomorphisms of free-abelian times free groups.
In Section \ref{inffixpersec} we count the orbits of the action of finite fixed (resp. periodic) points on infinite fixed (resp. periodic) points, proving that such conditions cannot be achieved in general for type I endomorphisms (in the sense of Delgado and Ventura) and that for type II endomorphisms $\Fix (\hat \varphi)$ (resp. $\Per(\hat\varphi)$)  has  $\sum\limits_{i=0}^m 2^i {m\choose i}$ $\Fix(\varphi)$-orbits (resp. $\Per(\varphi)$-orbits).  We also classify the infinite fixed points of automorphisms as attractors or repellers.
In Section \ref{dyninfpoints}  we study some dynamical aspects of infinite points of endomorphisms and automorphisms, proving that for type II endomorphisms every point is either periodic or wandering and that the dynamics of an automorphism is asymptotically periodic. Finally, in section \ref{furtherwk} we include some open questions.
\section{Uniform continuity of endomorphisms}
\label{unicontsec}
Since we are interested in the study of the continuous extension of endomorphisms to the completion, in this section we will obtain conditions for an endomorphism of a free-abelian group to be uniformly continuous. Whenever the metric is not mentioned, assume that it is the product metric given by taking the prefix metric in each direct factor. We present a proof of the following trivial lemma for sake of completeness.

\begin{lemma}
\label{produnicont}
Consider groups $G_i$ endowed with metrics $d_i$, for $i=1,2,3,4$. Let $\varphi_1:G_1\to G_3$ and $\varphi_2:G_2\to G_4$ be homomorphisms of groups. The homomorphism $\varphi:G_1\times G_2\to G_3\times G_4$ given by $(x,y)\varphi=(x\varphi_1,y\varphi_2)$ is uniformly continuous with respect to the product metrics $d$ and $d'$ if and only if $\varphi_1$ is uniformly continuous with respect to $d_1$ and $d_3$ and $\varphi_2$ is uniformly continuous with respect to $d_2$ and $d_4$.
\end{lemma}

\noindent \textit{Proof.}
Consider the homomorphism  $\varphi:G_1\times G_2\to G_3\times G_4$ given by $(x,y)\varphi=(x\varphi_1,y\varphi_2)$  and suppose that it is uniformly continuous with respect to the product metrics $d$ and $d'$s. Let $\varepsilon >0$ and take $\delta$ such that for every $(x_1,x_2),(y_1,y_2)\in G_1\times G_2$ such that $d((x_1,x_2),(y_1,y_2))<\delta$, we have that $d'((x_1\varphi_1,x_2\varphi_2),(y_1\varphi_1,y_2\varphi_2))<\varepsilon$. We know that for every $x_1,y_1\in G_1$ such that $d_1(x_1,y_1)<\delta$ we have $d_3(x_1\varphi_1,y_1\varphi_1)<\varepsilon,$ since $d((x_1,1),(y_1,1))=d_1(x_1,y_1)<\delta$ and so $$d_3(x_1\varphi_1,y_1\varphi_1)=d'((x_1\varphi_1,1\varphi_2),(y_1\varphi_1,1\varphi_2))<\varepsilon.$$

Conversely, if both $\varphi_1,\varphi_2$ are uniformly continuous, then taking $\varepsilon >0$, there are $\delta_i$ such that for every $x,y\in G_i$ such that $d_i(x,y)<\delta_i$, we have $d_{i+2}(x\varphi_i,y\varphi_i)<\varepsilon,$ for $i=1,2$. Taking $\delta=\min\delta_i$, we know that for every $(x_1,y_1),(x_2,y_2)\in G_1\times G_2$ such that $d((x_1,x_2),(y_1,y_2))<\delta$, then $d_1(x_1,y_1)<\delta \leq \delta_1$ and $d_2(x_2,y_2)<\delta \leq \delta_2$, thus  $d_{i+2}(x\varphi_i,y\varphi_i)<\varepsilon$ and  $$d'((x_1,x_2)\varphi,(y_1,y_2)\varphi)=d'((x_1\varphi_1,x_2\varphi_2),(y_1\varphi_1,y_2\varphi_2))<\varepsilon.$$
\qed\\

\begin{proposition}
\label{p1}
Let $u\in\mathbb Z^m$ and $\varphi:\mathbb Z^m \to \mathbb Z$ be a homomorphism given by $v\mapsto vu^T$. Then $\varphi$ is uniformly continuous if and only if $u$ has at most one nonzero entry.
\end{proposition}
\noindent \textit{Proof.}
If $u$=0, then $\varphi$ is  uniformly continuous.
It is clear that if $u$ has a single nonzero entry, then $\varphi$ is uniformly continuous.  Indeed, take $u\in \mathbb Z^m$ such that $u_k=\lambda\neq 0$ for some $k\in [m]$ and $u_j=0$ for all $j\in [m]\setminus \{k\}$.  Take $\varepsilon >0$. Set $\delta=\varepsilon $ and take $a,b\in \mathbb Z^m$ such that $d(a,b)<\delta$. Notice that $au^T=\lambda a_k$ and $bu^T=\lambda b_k$. If $a_k=b_k$, then $d(au^T,bu^T)=0<\varepsilon$. If not, since $d(a,b)<\delta$ then $d(a_i,b_i)<\delta$ for all $i\in [m]$. In particular $d(a_k,b_k)<\delta$. This means that $|a_k\wedge b_k|>log_2(\frac 1 \delta)$, i.e. $a_kb_k>0$ and $|a_k|,|b_k|>log_2(\frac 1\delta).$ But then, $au^Tbu^T=\lambda^2 a_kb_k>0$ and $|\lambda a_k|=|\lambda||a_k|\geq |a_k|>log_2(\frac 1 \delta)$. Similarly,  $|\lambda b_k|>log_2(\frac 1 \delta)$. This means that $d(au^T,bu^T)\leq d(a,b)<\delta=\varepsilon$

Suppose now that $u$ has at least one positive and one negative entry. Let $u_{i_1},\ldots,u_{i_r}$ be the nonnegative entries and $u_{j_1},\ldots,u_{j_s}$ be the negative entries and suppose w.l.o.g. that $\sum\limits_{x=1}^r u_{i_x}\geq \sum\limits_{x=1}^s u_{j_x}$. We will show that for every $\delta>0$, there are $v,w\in\mathbb Z^m$ such that $d(v,w)<\delta$ and $d(vu^T,wu^T)=1$ and so $\varphi$ is not uniformly continuous. Take $\delta>0$, $v$ such that
$$v_i=1+\lceil log_2(\frac 1\delta)\rceil \text{, for every $i\in [m]$}$$
and $w$ such that 
$$w_{i_k}=1+\left\lceil log_2\left(\frac 1\delta\right)\right\rceil \text{, for all $k\in[r]$}$$
and
$$ w_{j_k}=\sum\limits_{x=1}^r u_{i_x}\left(1+\left\lceil log_2\left(\frac 1\delta\right)\right\rceil\right)  \text{, for all $k\in[s]$}.$$
Then $d(v,w)<\delta$ since, for every $i\in [m]$, $v_iw_i>0$ and $|v_i|,|w_i|>log_2(\frac 1\delta)$ (notice that $\sum\limits_{x=1}^r u_{i_x}\geq 1$). Also, $vu^T=\sum v_iu_i\geq 0$, since we are assuming that $\sum\limits_{x=1}^r u_{i_x}\geq \sum\limits_{x=1}^s u_{j_x}$. We have that $$wu^T=\sum_{x=1}^r u_{i_x}\left(1+\left\lceil log_2\left(\frac 1\delta\right)\right\rceil\right)\left(1+\sum\limits_{x=1}^s u_{j_x}\right)\leq 0.$$ Thus, $d(vu^T,wu^T)=1.$

Now, suppose that $u\in(\mathbb Z_0^+)^m$ has at least two nonzero entries (the nonpositive case is analogous). Let $u_k$ be a minimal nonzero entry. As above, we will show that for every $\delta>0$, there are $v,w\in\mathbb Z^m$ such that $d(v,w)<\delta$ and $d(vu^T,wu^T)=1$ and so $\varphi$ is not uniformly continuous. Take $\delta>0$, $v$ such that  
$$v_i=1+\left\lceil log_2\left(\frac 1\delta\right)\right\rceil \text{ for every  $i \in [m]\setminus\{k\}$}\quad \text{and} \quad v_k=-1-\left\lceil log_2\left(\frac 1\delta\right)\right\rceil$$ 
and $w$ 
such that 
$$w_i=1+\left\lceil log_2\left(\frac 1\delta\right)\right\rceil \text{ for every  $i \in [m]\setminus\{k\}$}$$
and
$$w_k=-\sum\limits_{i\in[m]} u_i \left(1+\left\lceil log_2\left(\frac 1\delta\right)\right\rceil\right).$$
Now, $w_iv_i>0$ and $|w_i|,|v_i|>log_2(\frac 1 \delta)$ for every $i\in [m]$, so $d(v,w)<\delta.$ Also, 
$$vu^T=\left(1+\left\lceil log_2\left(\frac 1\delta\right)\right\rceil \right) \left(\sum\limits_{i\neq k} u_i -u_k\right)\geq 0$$ 
by minimality of $u_k$. We have that 
$$wu^T=\left(1+\left\lceil log_2\left(\frac 1\delta\right)\right\rceil \right) \left(\sum\limits_{i\neq k} u_i -u_k\sum u_i\right)\leq 0.$$ 
Thus, $d(vu^T,wu^T)=1.$
\qed\\

\begin{corollary}
\label{unicontfreeab}
Let $Q\in\mathbb \mc M_m(\mathbb Z)$ and $\varphi\in End(\mathbb Z^m)$ to be given by $u \mapsto uQ$. Then $\varphi$ is uniformly continuous if and only if every column $Q_i$ of $Q$ has at most one nonzero entry.
\end{corollary}
\noindent \textit{Proof.} Consider the homomorphisms $\varphi_i: \mathbb Z^m \to \mathbb Z$ defined by $u\mapsto uQ_i$. Then $\varphi(u)=(\varphi_1(u),\cdots, \varphi_m(u)).$ 
\qed\\

We are now capable of obtaining conditions of uniform continuity for endomorphisms of type I.
\begin{proposition}
\label{typeiuc}
Let $G=\mathbb Z^m\times F_n$, with $n>1$ and consider an endomorphism $\varphi$ of type I, mapping $(a,u)$ to $(aQ+\mathbf uP,u\Phi)$. Denote by $\psi$ the endomorphism of $\mathbb Z^m$ defined by $ a\mapsto {aQ}.$  Then the following conditions are equivalent:
\begin{enumerate}[i.]
\item $\varphi$ is uniformly continuous.
\item $P=0$, $\psi$ is uniformly continuous and $\Phi$ is either constant or injective.
\end{enumerate}
\end{proposition}
\noindent \textit{Proof.} $i.\Rightarrow ii.$  
Consider the alphabet of $F_n$ to be $\{x_1\ldots x_n\}$.  Suppose $P\neq 0$ and pick entries $p_{rs}\neq 0$ and $p_{ts}$ with $r\neq t$.  We will prove that $\varphi$ is not uniformly continuous, by showing that $\forall \delta>0$ there exists $X,Y\in G$ such that $d(X,Y)<\delta$ and $d(X\varphi,Y\varphi)=1$. We may assume $\delta\leq 1$, so pick such $\delta$ and, as usual,  set $q=1+\lceil log_2(\frac{1}{\delta})\rceil$. Take $\beta\in \mathbb Z$ to be such that $sgn(p_{ts}q)\neq sgn(p_{ts}q+\beta p_{rs})$ (if $p_{ts}=0$, put $\beta=1$). Let $X=(0,x_t^q)$  and $Y=(0,x_t^qx_r^\beta)$. To simplify notation, write $u$ and $v$ for the free parts of $X$ and $Y$, respectively, so $u=x_t^q$ and $v=x_t^qx_r^\beta$.

Since the free abelian parts coincide, $d(X,Y)=d(u,v)=2^{-q}<\delta$. 

We have that $d(( \mathbf uP,u\Phi),(vP,v\Phi))=\max\{d(u\Phi,v\Phi),d(\mathbf uP,vP)\}\geq d( \mathbf uP, vP)$. But 

\[ \mathbf uP =
\begin{bmatrix}
 0 & \cdots &   q   &  0 & \cdots & 0 
\end{bmatrix}
\begin{bmatrix}
    p_{11}  & \cdots & p_{1m}      \\
    \vdots \\
    p_{n1}     & \cdots & p_{nm}      
\end{bmatrix} 
=\begin{bmatrix}
   p_{ti}q
\end{bmatrix}_{i\in[m]}
\]
and 
\[\mathbf vP =
\begin{bmatrix}
0 & \cdots &    q & \cdots & \beta   & \cdots & 0 
\end{bmatrix}
\begin{bmatrix}
    p_{11}   & \cdots & p_{1m}      \\
    \vdots \\
    p_{n1}     & \cdots & p_{nm}      
\end{bmatrix} 
=\begin{bmatrix}
   p_{ti}q+\beta p_{ri}   
\end{bmatrix}_{i\in[m]}
\]
thus, $d(\mathbf{ u}P ,\mathbf{v}P)=\max \{d((\mathbf{ u}P)_i,(\mathbf{ v}P)_i)\}\geq d((\mathbf{ u}P)_s,(\mathbf{ v}P)_s)=d(  p_{ts}q ,p_{ts}q+\beta p_{rs})=1$,
by definition of $\beta$.  

The remaining conditions follow by Lemma \ref{produnicont}.

$ii.\Rightarrow i.$ This implication is obvious by Lemma \ref{produnicont} since both $a\mapsto aQ$ and $u\mapsto u\Phi$ are uniformly continuous.\qed\\

Now we deal with the type II endomorphisms. A reduced word $z=z_1\cdots z_n$, with $z_i \in A \cup A^{-1}$, is said to be \emph{cyclically reduced} if $z_1\neq z_n^{-1}$. Every word admits a decomposition of the form $z=w\tilde z w^{-1}$, where $\tilde z$ is cyclically reduced. The word $\tilde z$ is called the \emph{cyclically reduced core} of $z$. 

\begin{proposition}
\label{typeiiuc}
Let $G=\mathbb Z^m\times F_n$, with $n>1$ and consider an endomorphism $\varphi$ of type II, mapping $(a,u)$ to $(aQ+\mathbf uP,z^{{a\ell}^T+\mathbf{u}h^T})$. Denote by $\psi_1$ the endomorphism of $\mathbb Z^m$ defined by $ a\mapsto {aQ}$ and $\psi_2:\mathbb Z^m\to \mathbb Z$ the homomorphism defined by $ a\mapsto {a} \ell^T.$ Then the following conditions are equivalent:
\begin{enumerate}
\item $\varphi$ is uniformly continuous.
\item $\mathbf P=0$, $\mathbf h=0$ and both $\psi_1$ and $\psi_2$ are uniformly continuous.
\end{enumerate}
\end{proposition}
\noindent \textit{Proof.} $i.\Rightarrow ii.$ The proof that $\mathbf P=0$ is the same as in the previous proposition.

 Now, we will prove that if $ h\neq 0$, then for all $\delta >0$, there are $X,Y\in G$ such that $d(X,Y)<\delta$ and $d(X\varphi,Y\varphi)=1.$ Suppose then $ h\neq 0$ and pick entries $h_k\neq 0$ and $h_t$, with $t\neq k$ and some $\delta>0$. Set $q=1+\lceil log_2(\frac{1}{\delta})\rceil$ and take $X=(0,x_t^{q})$ and $Y=(0,x_t^{q}x_k^\alpha)$ for some $\alpha\in \mathbb Z$ such that $sgn(h_tq+\alpha h_k)\neq sgn(h_tq) $ (if $h_t=0$ put $\alpha=1$). Then $d(X,Y)<\delta$ and $1\geq d(X\varphi,Y\varphi)\geq d(z^{h_tq},z^{\alpha h_k+h_tq})=1$.
 
The proof that $\psi_1$ is uniformly continuous is analogous the the one in the previous proposition.

Now, suppose $\psi_2$ is not uniformly continuous. There exists $\varepsilon >0$ such that for every $\delta >0$, there are $a, b\in \mathbb Z^m$ such that $d( a,  b)<\delta$ and 
\begin{align}
\label{falhacont}
d\left(\sum\limits_{i\in [m]} a_i\ell_i,\sum\limits_{i\in[m]} b_i\ell_i\right)\geq\varepsilon.
\end{align}
 We now show that for every $\delta >0$, there are $X,Y\in G$ such that $d(X,Y)<\delta$ and $d(X\varphi,Y\varphi)\geq 2^{-|w|-|\tilde z|{\lceil log_2(\frac 1\varepsilon)\rceil}}$, where $\tilde z$ is the cyclically reduced core of $z$ and $w$ is such that $z=w\tilde z w^{-1}$. Notice that $\tilde z$ and $w^{-1}$ don't share a prefix. Take $\delta>0$ and take $ a, b\in \mathbb Z^m$ such that $d(a,b)<\delta$ satisfying (\ref{falhacont}). Now, consider $X=(a,1)$ and $Y=(b,1)$. Clearly $d(X,Y)=d( a, b)<\delta$ and 
 $$d(X\varphi,Y\varphi)=d\left(\left(aQ,z^{a\ell^T}\right),\left(bQ,z^{b\ell^T}\right)\right)\geq d\left(z^{{a\ell^T}},z^{{b\ell^T}}\right)$$
  We know that (\ref{falhacont}) holds, so, either $$a\ell^Tb\ell^T=\sum\limits_{i\in[m]} a_i\ell_i\sum\limits_{i\in[m]} b_i\ell_i\leq 0$$ and in that case $d\left(z^{a\ell^T},z^{b\ell^T}\right)=1$, or
   $$\sum\limits_{i\in[m]} a_i\ell_i\sum\limits_{i\in[m]} b_i\ell_i> 0 \text{ and }  2^{-\min \left\{\left|a\ell^T\right|,\left|b\ell^T\right|\right\}}\geq\varepsilon,$$
    which means that
     $$\min \left\{\left|a\ell^T\right|,\left|b\ell^T\right|\right\}\leq log_2 \left(\frac 1\varepsilon\right).$$ 
In this case, we have that
 \begin{align*}
 d(z^{a\ell^T},z^{b\ell^T})&=2^{-|z^{a\ell^T}\wedge z^{b\ell^T}|}=2^{-|w\tilde z^{a\ell^T}w^{-1}\wedge w\tilde z^{b\ell^T}w^{-1} |}\\
 &=2^{-|w|-|\tilde z|{\min\{|a\ell^T|,|b\ell^T|\}}}\geq 2^{-|w|-|\tilde z|{\lceil log_2(\frac 1\varepsilon)\rceil}}.
\end{align*}

$ii.\Rightarrow i.$ Straightforward.\qed\\

So,  type I uniformly continuous endomorphisms are of the form $(a,u)\mapsto (aQ,u\phi)$ where $Q\in \mathcal M_m(\mathbb Z)$ has at most one nonzero entry in each column and $\phi\in End(F_n)$ is either constant or injective. Type II endomorphisms are uniformly continuous if they map $(a,u)$ to $(aQ,z^{\lambda a_k})$ where $Q\in \mathcal M_m(\mathbb Z)$ has at most one nonzero entry in each column, $k\in [m]$, $1\neq z\in F_n$ is not a proper power and $0\neq\lambda\in\mathbb Z$. Notice that $\lambda \neq 0$, since by definition of a type II endomorphism we have that $\ell\neq 0.$

\begin{remark}
\label{ucautos}
In \cite{[DV13]}, the authors prove that an endomorphism $\varphi\in End(\mathbb Z^m\times F_n)$ is an automorphism if it is of type I and of the form $(a,u) \mapsto (aQ+\mathbf u P,u\phi)$, with $\phi\in Aut(F_n)$ and $ Q\in GL_m(\mathbb Z)$. In the case where $\varphi$ is uniformly continuous, then every entry of $Q$ is either $1$ or $-1$ and $Q=AD$, where $D$ is diagonal and $A$ is a permutation matrix. There are $2^mm!$ such matrices, which we call \emph{uniform}. So, a uniformly continuous automorphism of $\mathbb Z^m\times F_n$ is defined as $(a,u)\mapsto (aQ,u\phi)$, where $Q$ is a uniform matrix and $\phi\in Aut(F_n)$. 
\end{remark}


\section{Coarse-median preserving endomorphisms}
\label{cmpendos}

A metric space $(X,d)$ is said to be a \emph{median space} if, for all $x,y,z\in X$, there is some 
unique point $\mu(x,y,z)\in X$, known as the \emph{median} of $x,y,z$, 
such that $d(x,y)=d(x,\mu(x,y,z))+d(\mu(x,y,z),y);$ $d(y,z)=d(y,\mu(x,y,z))+d(\mu(x,y,z),z);$ and $d(z,x)=d(z,\mu(x,y,z))+d(\mu(x,y,z),x)$. We call $\mu:X^3\to X$ the \emph{median operator} of the median space $X$. 

Coarse median spaces were introduced by Bowditch in \cite{[Bow13]}. Following the equivalent definition given in \cite{[NWZ19]}, we say that, given a metric space $X$, a \emph{coarse median on $X$} is a ternary operation $\mu:X^3\to X$ satisfying the following: 

there exists a constant $C\geq 0$ such that, for all $a,b,c,x\in X$, we have that 
\begin{enumerate}
\item $\mu(a,a,b)=a \text{ and } \mu(a,b,c)=\mu(b,c,a)=\mu(b,a,c);$
\item $d(\mu(\mu(a,x,b),x,c),\mu(a,x,\mu(b,x,c)))\leq C;$
\item $d(\mu(a,b,c),\mu(x,b,c))\leq Cd(a,x)+C.$
\end{enumerate} 

Given a group $G$, a \emph{word metric} on $G$ measures the distance of the shortest path in the Cayley graph of $G$ with respect to some set of generators, i.e., for two elements $g,h\in G$, we have that $d(g,h)$ is the length of the shortest word whose letters come from the generating set representing $g^{-1}h$.
Following the definitions in \cite{[Fio21]}, two coarse medians $\mu_1,\mu_2:X^3\to X$ are said to be at \emph{bounded distance} if there exists some constant $C$ such that $d(\mu_1(x,y,z),\mu_2(x,y,z))\leq C$ for all $x,y,z\in X$, and a \emph{coarse median structure} on $X$ is an equivalence class $[\mu]$ of coarse medians pairwise at bounded distance. When $X$ is a metric space and $[\mu]$ is a coarse median structure on $X$, we say that $(X,[\mu])$ is a \emph{coarse median space.} Following Fioravanti's definition, a \emph{coarse median group} is a pair $(G,[\mu])$, where $G$ is a finitely generated group with a word metric $d$ and $[\mu]$ is a $G$-invariant coarse median structure on $G$, meaning that for each $g\in G$, there is a constant $C(g)$ such that $d(g\mu(g_1,g_2,g_3),\mu(gg_1,gg_2,gg_3)\leq C(g)$, for all $g_1,g_2,g_3\in G$. The author in \cite{[Fio21]} also remarks that this definition is stronger than the original definition from \cite{[Bow13]}, that did not require $G$-invariance. Despite being better suited for this work, it is not quasi-isometry-invariant nor commensurability-invariant, unlike Bowditch's version.

In $\mathbb Z^m$, the $\ell_1$ metric is defined by $$d_{\ell_1}(a,b)=\sum_{i=1}^m|a_i-b_i|.$$Consider $\mathbb Z^m$ endowed with the $\ell_1$ metric. Then, given three points $a,b,c\in \mathbb Z^m$, consider $\mu(a,b,c)$ to be the point having in the component $i$ the (numerical) median of $\{a_i,b_i,c_i\}$. Then $\mu$ is a median operator.

Now, a $K$-hyperbolic group is such that there is some $K\geq 0$ for which every geodesic triangle has a $K$-center, i.e., a point that, up to a bounded distance, depends only on the vertices, and is $K$-close to every edge of the triangle. Given three points, the operator that associates the three points to the $K$-center of a geodesic triangle they define is coarse median. In fact, by Theorem 4.2. in \cite{[NWZ19]} it is the only coarse-median structure that we can endow $X$ with.

Given two groups $G_1$ and $G_2$ endowed with coarse median operators $\mu_1$ and $\mu_2$, then it is easy to check that the operator $\mu:(G_1\times G_2)^3\to G_1\times G_2$ defined by $\mu((x_1,x_2),(y_1,y_2),(z_1,z_2))=(\mu_1(x_1,y_1,z_1),\mu_2(x_2,y_2,z_2))$ is coarse median.This can naturally be extended to direct products of $n$ factors. We refer to $\mu$ as the \emph{product coarse median operator}.\\

Given a coarse median group $(G,[\mu])$, an automorphism $\varphi\in \Aut(G)$ is said to be \emph{coarse-median preserving} if it fixes $[\mu]$, i.e., when there is some constant $C\geq 0$ such that for all $g_i'$s, we have that $$d(\mu(g_1,g_2,g_3)\varphi,\mu(g_1\varphi,g_2\varphi,g_3\varphi))\leq C,$$ with respect to some word metric $d$. This can naturally be defined for general endomorphisms, not necessarily bijective.\\

\begin{remark}
\label{cmpfiora}
Example 2.26 in \cite{[Fio21]} says that, in $\mathbb Z^m$, when taking $\mu$ to be the median operator associated with the $\ell_1$ metric, the coarse median preserving automorphisms are the ones given by uniform matrices, which correspond to the uniformly continuous automorphisms of $\mathbb Z^m$ when the product metric is taken with the prefix metric in each component (recall Remark \ref{ucautos}). In the case of hyperbolic groups, with the coarse median defined above, every automorphism is coarse-median preserving. In fact, it is proved in \cite{[Car21]} that given a hyperbolic group $G$  and an endomorphism $\varphi\in \End(G)$, then the bounded reduction property (BRP) holds for $\varphi$ if and only if $\varphi$ is coarse-median preserving.
\end{remark}

We will now see that, in some sense, coarse-median preserving endomorphisms coincide with the uniformly continuous ones for $\mathbb Z^m\times F_n.$ We start with a fairly obvious technical lemma. Observe that, given two groups $G_1=\langle A\rangle$ and $G_2=\langle B\rangle$ endowed with word metrics $d_1$ and $d_2$, respectively, then the product metric $d$ is a word metric for the generators $A^1\times B^1$.
\begin{lemma}
\label{cmpprod}
Let $G_i$ be groups endowed with  word metrics $d_i$ and coarse medians $\mu_i$, and consider endomorphisms $\phi_i\in \End(G_i),$  for $i=1,\ldots, k$. The endomorphism $\varphi\in \End(G_1\times \cdots\times G_k)$ defined by $(x_1,\ldots x_k)\mapsto (x_1\phi_1,\ldots, x_k\phi_k)$ is coarse-median preserving with respect to the product coarse median operator and the product metric $d$ if and only if $\phi_i$ is coarse-median preserving for $\mu_i$ and $d_i$ for every $i\in [k]$.
\end{lemma}
\noindent \textit{Proof.}	Suppose that, for every $i\in[k]$, the endomorphism $\phi_i$ is coarse-median preserving and take $C=\max\{C_i\mid i \in [k]\}$, where $C_i$ is the constant given by this property for $\phi_i$. Then, we have that, for every $i\in [k]$,
\begin{align}
\label{eachcmp}
d_i(\mu_i(x_i,y_i,z_i)\phi_i,\mu_i(x_i\phi_i,y_i\phi_i,z_i\phi_i)\leq C_i\leq C,
\end{align}
for all $x_i,y_i,z_i\in G_i.$
Thus, for all $(x_1,\ldots,x_k),(y_1,\ldots,y_k), (z_1,\ldots, z_k)\in G_1 \times \cdots\times G_k$, we have that  
\begin{align*}
&\mu((x_1,\ldots,x_k)\varphi,(y_1,\ldots,y_k)\varphi, (z_1,\ldots, z_k)\varphi)\\
=&(\mu_1(x_1\phi_1,y_1\phi_1,z_1\phi_1),\ldots, \mu_k(x_k\phi_k,y_k\phi_k,z_k\phi_k))
\end{align*}
and 
$$\left(\mu((x_1,\ldots x_k),(y_1,\ldots, y_k),(z_1,\ldots,z_k))\right)\varphi=((\mu_1(x_1,y_1,z_1))\phi_1,\ldots, (\mu_k(x_k,y_k,z_k))\phi_k).$$ 
From (\ref{eachcmp}), it follows that 
$
\mu((x_1,\ldots,x_k)\varphi,(y_1,\ldots,y_k)\varphi, (z_1,\ldots, z_k)\varphi)$ is $C$-close to $\left(\mu((x_1,\ldots x_k),(y_1,\ldots, y_k),(z_1,\ldots,z_k))\right)\varphi)$ with respesct to the product metric $d$,
and so $\varphi$ is coarse-median preserving.

The converse is proved similarly.
\qed\\

\begin{lemma}
\label{facmpuc}
An endomorphism $\varphi\in \End(\mathbb Z^m)$ is uniformly continuous with respect to the product metric $d$ obtained by taking the prefix metric in each direct factor if and only if it is coarse-median preserving for the median operator $\mu$ induced by the metric $\ell_1$ in $\mathbb Z^m$.
\end{lemma}
\noindent \textit{Proof.}
Let $\varphi\in \End(\mathbb Z^m)$ be given by $Q\in \mc M_m(\mathbb Z)$ and suppose that it is coarse-median preserving. We proceed in a similar way to what we did in the proof of Proposition \ref{p1}. If $\varphi$ is not uniformly continuous with respect to $d$, then there is some column $Q_j$ and nonzero entries $q_{rj}$ and $q_{sj}$ with $s>r$. If they are both positive, suppose w.l.o.g. that $q_{sj}>q_{rj}$. For every $n\in \N$, let $x^{(n)}\in \mathbb Z^{m}$ such that $x^{(n)}_{r}=1$, $x^{(n)}_{s}=n$ and all the other entries are $0$; $y^{(n)}\in \mathbb Z^{m}$ such that $y^{(n)}_{r}=1$ and all the other entries are $0$; and $z^{(n)}\in \mathbb Z^{m}$ such that $z^{(n)}_{r}=1+2nq_{sj}$, $z^{(n)}_{s}=-1$ and all the other entries are $0$. Then 
$\mu(x^{(n)},y^{(n)},z^{(n)})$ has $1$ in the $r$-th entry and all the other entries are $0$, so $[\mu(x^{(n)},y^{(n)},z^{(n)})\varphi]_j=q_{r_j}$. But $$[\mu(x^{(n)}\varphi,y^{(n)}\varphi,z^{(n)}\varphi)]_j=\mu(q_{rj}+nq_{sj},q_{rj},q_{rj}+(2nq_{rj}-1)q_{sj})=q_{rj}+nq_{sj}$$
and so $d_{\ell^1}(\mu(x^{(n)}\varphi,y^{(n)}\varphi,z^{(n)}\varphi),\mu(x^{(n)},y^{(n)},z^{(n)})\varphi)\geq nq_{sj}$, and that contradicts the fact that $\varphi$ is coarse-median preserving. If both $q_{rj}$ and $q_{sj}$ are negative, we can reach a contradiction in the same way. If, suppose $q_{rj}>0$ and $q_{sj}>0$, putting, for all $n\in \mathbb N$, $x^{(n)}\in \mathbb Z^{m}$ such that 
$x^{(n)}_{r}=-nq_{sj}$, $x^{(n)}_{s}=nq_{rj}$ and all the other entries are $0$; 
$y^{(n)}\in \mathbb Z^{m}$ such that $y^{(n)}_{r}=-2nq_{sj}$, $y^{(n)}_{s}=-nq_{rj}$ and all the other entries are $0$; 
and 
$z^{(n)}\in \mathbb Z^{m}$ such that $z^{(n)}_{r}=-3nq_{sj}$, $z^{(n)}_{s}=4nq_{rj}$ and all the other entries are $0$, we get that $
\mu(x^{(n)},y^{(n)},z^{(n)})$ has $-2nq_{sj}$ in the $r$-th entry, $nq_{rj}$ in the $s$-th entry and all the other entries are $0$, so $[\mu(x^{(n)},y^{(n)},z^{(n)})\varphi]_j=-nq_{rj}q_{sj}$. 
But $$[\mu(x^{(n)}\varphi,y^{(n)}\varphi,z^{(n)}\varphi)]_j=\mu(0,-3nq_{rj}q_{sj},nq_{rj}q_{sj})=0$$
and so $d_{\ell^1}(\mu(x^{(n)}\varphi,y^{(n)}\varphi,z^{(n)}\varphi),\mu(x^{(n)},y^{(n)},z^{(n)})\varphi)\geq nq_{rj}q_{sj}$, and that contradicts the fact that $\varphi$ is coarse-median preserving. 

Hence, we have that $\varphi$ is uniformly continuous with respect to $d$.

Conversely, suppose that $\varphi$ is uniformly continuous with respect to $d$ and recall the notation introduced in Section \ref{inffixpersec}. Consider three points $a,b,c\in \mathbb Z^m$. Then $$\mu(aQ,bQ,cQ)=[\mu(\lambda_ia_{\alpha_i},\lambda_ib_{\alpha_i},\lambda_ic_{\alpha_i})]_{i\in [m]}=[\lambda_i\mu(a_{\alpha_i},b_{\alpha_i},c_{\alpha_i})]_{i\in [m]}=\mu(a,b,c)Q.$$
Thus $\varphi$ is coarse-median preserving.
\qed

We are now able to prove the main result of this section.
\begin{theorem}
\label{cmpfatf}
An endomorphism $\varphi\in \End(\mathbb Z^m\times F_n)$ is uniformly continuous with respect to the product metric $d$ obtained by taking the prefix metric in each direct factor if and only if it is coarse-median preserving for the product coarse median $\mu$ obtained by taking the median operator $\mu_1$ induced by the metric $\ell_1$ in $\mathbb Z^m$ and the coarse median operator $\mu_2$ given by hyperbolicity in $F_n$.
\end{theorem}

\noindent \textit{Proof.}	
An endomorphism of $\mathbb Z^m \times F_n$ is of the form $(a,u)\mapsto (aQ+\mathbf uP, (a,u)\psi)$, where $Q\in \mc M_m(\mathbb Z)$, $P\in \mc M_{n\times m}(\mathbb Z)$ and $\psi$ is a homomorphism from $\mathbb Z^m\times F_n$ to $F_n$. We start by proving that if $\varphi$ is coarse-median preserving for $\mu$, then, $P=0$. 
Suppose then that $P\neq 0$ and pick entries $p_{rs}\neq 0$ and $p_{ts}$ with $r\neq t$. Take $\{x_1,\ldots,x_n\}$ to be a basis of $F_n$ and let $d_2$ be the word metric defined by this set of generators. Now, we have that  
$$(\mu((0,x_t), (0,x_tx_r^{C+1}),(0,x_r^{C+1}x_t)))\varphi=(0,x_t)\varphi=([p_{ti}]_{i\in[m]},(0,x_t)\psi)$$
and 
\begin{align*}
&\mu\left((\mathbf{x_t}P,(0,x_t)\psi),(\mathbf{x_tx_r}^{C+1}P,(0,x_tx_r^{C+1})\psi),(\mathbf{x_r}^{C+1}\mathbf{x_t}P,(0,x_r^{C+1}x_t)\psi)\right)\\
=&\left(\mathbf{x_tx_r}^{C+1}P,\mu_2((0,x_t)\psi,(0,x_tx_r^{C+1})\psi,(0,x_r^{C+1}x_t)\psi)\right)\\
=&\left([p_{ti}+p_{ri}(C+1)]_{i\in[m]},\mu_2((0,x_t)\psi,(0,x_tx_r^{C+1})\psi,(0,x_r^{C+1}x_t)\psi)\right).
\end{align*}

But taking the $\ell_1$ metric in $\mathbb Z^m$, we have that $$d_{\ell_1}([p_{ti}]_{i\in[m]},[p_{ti}+p_{ri}(C+1)]_{i\in[m]})\geq |(C+1)p_{rs}|>C,$$ which contradicts the fact that $\varphi$ is coarse-median preserving with constant $C$. So, $P$ must be zero.

If $\varphi$ is a type I coarse-median preserving endomorphism for $\mu$, then $\varphi$ is of the form $(a,u)\mapsto (aQ, u\phi)$, where $Q\in \mc M_m(\mathbb Z)$ and $\phi\in \End(F_n).$ By Lemma \ref{cmpprod}, we get that both the endomorphism $\theta$ of $\mathbb Z^m$ defined by $Q$ and $\phi\in \End(F_n)$ must be coarse-median preserving. By Lemma \ref{facmpuc}, $\theta$ is uniformly continuous and by Remark \ref{cmpfiora}, we have that the BRP holds for $\phi$, and so $\phi$ is uniformly continuous. By Proposition \ref{typeiuc}, $\varphi$ is uniformly continuous for the metric $d$. 

Conversely, if $\varphi$ is a type I uniformly continuous with respect to the metric $d$, then, by Proposition  \ref{typeiuc}, $\varphi$ is of the form $(a,u)\mapsto (aQ, u\phi)$, where  $Q\in \mc M_m(\mathbb Z)$ has at most one nonzero entry in each column and $\phi\in \End(F_n)$ is uniformly continuous for the prefix metric. Combining Lemmas \ref{cmpprod} and \ref{facmpuc} with Remark \ref{cmpfiora}, we get that $\varphi$ is coarse-median preserving for $\mu$. 

Now we deal with the type II endomorphisms. Suppose that $\varphi$ is a coarse-median preserving type II endomorphism with constant $C\geq 0$. Then, it must be of the form $(a,u)\mapsto (aQ, z^{a\ell^T+\textbf uh^T})$, for some $\mathbf 0\neq \ell\in \mathbb Z^m$ and $h\in \mathbb Z^n$. We will prove, proceeding in the same way we did above, that $h$ must be zero and that $\ell$ must have at most one nonzero entry. Suppose that $h$ is not zero and pick $h_s\neq 0$ and $h_t$, with $t\neq s$. Then 
$$(\mu((0,x_t), (0,x_tx_r^{C+1}),(0,x_r^{C+1}x_t)))\varphi=(0,x_t)\varphi=(0,z^{h_t})$$
and 
\begin{align*}
&\mu\left((0,z^{h_t}),(0,z^{h_t+(C+1)h_s}),(0,z^{h_t+(C+1)h_s})\right)\\
=&(0,z^{h_t+(C+1)h_s})
\end{align*}

But we have that $$d_2(z^{h_t},z^{h_t+(C+1)h_s})>C,$$ which contradicts the fact that $\varphi$ is coarse-median preserving with constant $C$. So, $h$ must be zero.
To prove that $\ell$ must have at most one nonzero entry, we proceed as in the proof of Lemma \ref{facmpuc}. For all $n\in \N$, we define $x^{(n)}$, $y^{(n)}$ and $z^{(n)}$ in the same way replacing $q_{rj}$ and $q_{sj}$ by $l_r$ and $l_s$, respectively. We have that the distance of the free component between $\mu((x^{(n)},1),(y^{(n)},1),(z^{(n)},1))\varphi$ and $\mu((x^{(n)},1)\varphi,(y^{(n)},1)\varphi,(z^{(n)},1)\varphi)$ will be unbounded, contradicting coarse-median preservation of $\varphi$. 

Conversely, if $\varphi$ is a uniformly continuous type II endomorphism, then it is of the form $(a,u)\mapsto(aQ,z^{\lambda a_k})$, for some $Q\in \mc M_m(\mathbb Z)$ with at most on nonzero entry in each column, $k\in[m]$, $1\neq z\in F_n$, which is not a proper power and $0\neq \lambda\in \mathbb Z$. Then, taking $(a,u),(b,v), (c,w)\in \mathbb Z^m\times F_n$, and letting $\mu_3:\mathbb Z^3\to \mathbb Z$ denote the usual median in $\mathbb Z$, we have that 
$$\mu((a,u),(b,v), (c,w))\varphi=(\mu_1(a,b,c)Q, z^{\lambda\mu_3(a_k,b_k,c_k)}) $$
and 
$$\mu((a,u)\varphi,(b,v)\varphi, (c,w)\varphi)=(\mu_1(aQ,bQ,cQ), \mu_2(z^{\lambda a_k},z^{\lambda b_k},z^{\lambda c_k})).$$
Since $Q$ defines a uniformly continuous endomorphism of $\mathbb Z^m$, it follows from Lemma \ref{facmpuc} that it is coarse-median preserving, thus $\mu_1(aQ,bQ,cQ)$ and $\mu_1(a,b,c)Q$ are close. 
Letting $z=w\tilde z w^{-1}$, where $\tilde z$ is the cyclically reduced core of $z$, we have that $w\tilde z^{\mu_3(\lambda a_k,\lambda b_k,\lambda c_k)}=w\tilde z^{\lambda\mu_3( a_k, b_k, c_k)}$ belongs to every edge of the geodesic triangle defined by $z^{\lambda a_k},z^{\lambda b_k}$ and $z^{\lambda c_k}$. Since $$d_2(z^{\lambda\mu_3(a_k,b_k,c_k)},w\tilde z^{\lambda\mu_3( a_k, b_k, c_k)})=d_2(w\tilde z^{\lambda\mu_3(a_k,b_k,c_k)}w^{-1},w\tilde z^{\lambda\mu_3( a_k, b_k, c_k)})=|w|,$$ we have that the distance between $z^{\lambda\mu_3(a_k,b_k,c_k)}$ and $ \mu_2(z^{\lambda a_k},z^{\lambda b_k},z^{\lambda c_k})$ is bounded.
\qed

\section{Infinite fixed and periodic points}
\label{inffixpersec}

In this section, we will study \emph{infinite} fixed and periodic points, i.e., considering a uniformly continuous endomorphism $\varphi$ and its continuous extension to the completion, $\hat\varphi$, we study the behaviour of fixed and periodic points of $\hat\varphi$ belonging to the boundary of the group. An infinite fixed point is said to be \emph{singular} if it belongs to the topological closure $(\Fix(\varphi))^c$ of $\Fix(\varphi)$ and \emph{regular} if it doesn't. We denote by $\Sing(\hat\varphi)$ (resp. $\Reg(\hat\varphi)$) the set of all singular (resp. regular) infinite fixed points of $\hat\varphi.$

We start by obtaining finiteness conditions on infinite fixed and periodic points, counting the number of $\Fix(\varphi)$-orbits of $\Fix(\hat\varphi)$ and of $\Per(\varphi)$-orbits of $\Per(\hat\varphi)$, under natural actions of $\Fix(\varphi)$ and $\Per(\varphi)$ on $\Fix(\hat\varphi)$ and $\Per(\hat\varphi)$, respectivelym when $\varphi$ is a uniformly continuous type II endomorphism. In particular, the number of orbits is finite, which can be seen as some sort of infinite version of Proposition 6.2 in \cite{[DV13]}.

Finally, we classify infinite fixed points as attractors or repellers.

\subsection{Finiteness Conditions on Infinite Fixed Points}
\label{finitinffixed}
Let $\varphi$ be a type I uniformly continuous endomorphism of $\mathbb Z^m \times F_n$, with $n>1$. Then $\varphi:\mathbb Z^m \times F_n \to \mathbb Z^m\times F_n$ is given by $ (a,u)\mapsto (aQ,u\phi)$, for some $ Q\in \mathcal M_m(\mathbb Z)$ such that every column of $ Q$ contains at most one nonzero entry and some either constant or injective $\phi\in End(F_n)$. Consider $\varphi_1\in End(\mathbb Z^m)$  to be defined as $a\mapsto aQ$. Clearly, $\Fix(\varphi)=\Fix(\varphi_1)\times \Fix(\phi)$ and it is finitely generated. 

Let $\hat\varphi:\widehat{\mathbb Z^m}\times \widehat{F_n}\to \widehat{\mathbb Z^m}\times \widehat{F_n}$ be its continuous extension to the completion. By uniqueness of the extension, we have that $\hat \varphi$ is given by $(a,u)\mapsto (a\widehat\varphi_1,u\widehat \phi)$.  Then, we have that  $\text{Fix}(\hat\varphi)= \text{Fix}(\hat \varphi_1)\times \text{Fix}(\hat \phi)$; $\text{Sing}(\hat\varphi)=\text{Sing}(\hat \varphi_1)\times \text{Sing}(\hat\phi)$ and $\text{Reg}(\hat\varphi)=\text{Fix}(\hat \varphi_1)\times \text{Reg}(\hat \phi)\cup \text{Reg}(\widehat \varphi_1)\times\text{Fix}(\hat\phi)$. There is no hope of finding a finiteness condition in this case that holds in general since, if $n=2$ and $\phi$ is the identity mapping, which is injective, then $\text{Sing}(\hat \phi)$ is uncountable, thus so are both $\text{Reg}(\hat\varphi)$ and $\text{Sing}(\hat\varphi)$.\\

However, we will see that this is not the case when dealing with type II endomorphisms.

Let $\varphi\in End(\mathbb Z^m)$ defined by $a\mapsto aQ$ be a uniformly continuous endomorphism. From Corollary \ref{unicontfreeab} we know that $Q$ has at most one nonzero entry in each column. Given a column $Q_j$, if $Q_j\neq 0$, we call $\lambda_j$ its nonzero entry and denote by $\alpha_j$ the corresponding row, so that $q_{ij}=\lambda_i$ if $i=\alpha_j$ and $q_{ij}=0$ otherwise. If $Q_j=0$ we put $\lambda_j=0$ and $\alpha_j=1$. Then, we have that $[a_i]_{i\in [m]}$ is mapped to $[\lambda_ia_{\alpha_i}]_{i\in [m]}$ and $a\in \Fix(\varphi)$ if $a_i=\lambda_ia_{\alpha_i}$ for every $i\in [m].$ Take $\hat\varphi_1:\widehat{\mathbb Z^m}\to \widehat{\mathbb Z^m}$ to be the continuous extension of $\varphi_1$ to the completion.

Given $a\in \widehat{\mathbb Z^m}$, we set $i_1,\ldots,i_r$ to be the indices such that $a_{i_1}=\ldots =a_{i_r}=+\infty$;  $j_1,\ldots,j_s$ to be the indices such that $a_{j_1}=\ldots= a_{j_s}=-\infty$ and $k_1,\ldots,k_t$ to be the indices such that $a_{k_1},\ldots a_{k_t}\not\in\{-\infty,+\infty\}$. Define, for every $n\in \mathbb N$,  $a_n\in \mathbb Z^m$ such that 
$$a_{n_{i_l}}=n \text{, for $l$ in $[r]$;} \quad a_{n_{j_l}}=-n \text{, for $l$ in $[s]$} \quad  \text{ and }\quad a_{n_{k_l}}=a_{k_l} \text{ for $l$ in $[t]$.} $$
We have that, given $\varepsilon >0$, for every $n>log_2(\frac 1\varepsilon)$, $d(a_n,a)<\varepsilon$, so $(a_n)\to a.$ Thus, $a\hat\varphi=(\lim a_n)\hat\varphi=\lim (a_n\varphi)$. Since $(a_n\varphi)$ is such that $(a_n\varphi)_i=\lambda_ia_{\alpha_i}$ if $a_{\alpha_i}\not\in\{+\infty,-\infty\}$, $(a_n\varphi)_i=n\lambda_i$ if $a_{\alpha_i}=+\infty$, and $(a_n\varphi)_i=-n\lambda_i$ if $a_{\alpha_i}=-\infty$, we have that $a\hat\varphi=[\lambda_ia_{\alpha_i}]$, assuming that $0\times\infty=0$.
So, $a\in \Fix(\hat\varphi)$ if and only if $a_i=\lambda_i a_{\alpha_i}$ for every $i\in[m]$.\\

Defining the sum of an integer with infinite in the natural way, we have that the subgroup $\Fix(\varphi)\leq \mathbb Z^m\times F_n$ acts on $\Fix(\hat\varphi)$ by left multiplication. Given $a\in \Fix(\varphi)$ and $b\in \Fix(\hat\varphi)$, then $(a+b)\hat\varphi=a\varphi+b\hat\varphi=a+b\in \Fix(\hat\varphi)$. We now count the orbit of this action.

\begin{proposition}
\label{finitefa}
Let $\varphi\in \End(\mathbb Z^m)$ be a uniformly continuous endomorphism. Then $\Fix (\hat\varphi)$ has $\sum\limits_{i=0}^m 2^i {m\choose i}$ $\Fix(\varphi)$-orbits.
\end{proposition}
\noindent \textit{Proof.}
Let $a\in \Fix(\hat\varphi)$ and define $r,s,t\in \{0,\ldots, m\}$, and $i_l,j_l,k_l$ as above. We will prove that, for $b\in \Fix(\hat\varphi)$, we have that  $b\in (\Fix \varphi)a$ if and only if $b_{i_l}=a_{i_l}$ for every $l\in [r]$, $b_{j_l}=a_{j_l}$, for every $l\in [s]$ and $b_{k_l}\not\in \{+\infty,-\infty\}$, for every $l\in [t]$, i.e., if their infinite entries coincide. If that is the case, then every orbit is defined by the position and the signal of their infinite entries. Obviously, for $i\in \{0,\ldots,m\}$, there are $m\choose i$ choices for $i$ infinite entries and each of them can be $+\infty$ or $-\infty$, hence the $2^i$ factor.

Start by supposing that $b\in \Fix(\hat\varphi)$ is such that $b\in (\Fix\varphi)a$. Then, there is some $c\in \Fix\varphi$ such that $b=c+a$. This means that for every $l\in [r]$, we have that $b_{i_l}=c_{i_l}+a_{i_l}=c_{i_l}+ (+\infty)=+\infty$, since $c\in \mathbb Z^m$. Similarly, we have that $b_{j_l}=a_{j_l}$, for every $j\in[s]$ and $b_{k_l}\not\in \{+\infty,-\infty\}.$ It is clear that $|b_{k_l}|<\infty$ for $l\in [t].$

Now, suppose $b\in \Fix(\hat\varphi)$ is such that $b_{i_l}=+\infty$, for $l\in[r]$, $b_{j_l}=-\infty$, for $l\in [s]$ and $b_{k_l}\not\in\{+\infty,-\infty\}$, for $l\in [t]$. Consider $c\in \mathbb Z^m$ defined by $c_{k_l}=b_{k_l}-a_{k_l}$ and all other entries are $0$. Clearly $b=c+a$. We only have to check that $c\in \Fix(\varphi),$ i.e., $\lambda_ic_{\alpha_i}=c_i$ for every $i\in [m]$. For $i$ such that $a_i=\pm \infty$, we have that $c_i=0$ and $a_{\alpha_i}=sgn(\lambda_i)a_i=\pm\infty$, which implies that $c_{\alpha_i}=0$. If not, $c_{i}=b_i-a_i=\lambda_i(b_{\alpha_i}-a_{\alpha_i})=\lambda_ic_{\alpha_i}$ and we are done. 
\qed\\

Now, let $\varphi$ be a type II uniformly continuous endomorphism of $\mathbb Z^m \times F_n$, with $n>1$. Then $\varphi:\mathbb Z^m \times F_n \to \mathbb Z^m\times F_n$ is given by $(a,u)\mapsto (aQ,z^{\lambda a_k})$, for some $ Q\in \mathcal M_m(\mathbb Z)$ such that every column of $ Q$ contains at most one nonzero entry, $0\neq\lambda\in \mathbb Z$ and $k\in [m]$. Consider $\varphi_1\in End(\mathbb Z^m)$  to be defined as $a\mapsto aQ$ and $\varphi_2: \mathbb Z^m\to F_n$ that maps $a$ to $z^{\lambda a_k}$, which are both uniformly continuous. Observe that $\Fix (\varphi)=\{(a,a\varphi_2)\mid a\in \Fix(\varphi_1)\}$ and it is finitely generated (see Proposition 6.2 in \cite{[DV13]}).

By uniqueness of extension, we have that $\hat\varphi: \widehat{\mathbb Z^m}\times \widehat{F_n}\to \widehat{\mathbb Z^m}\times \widehat{F_n}$ is defined by $(a,u)\mapsto (a\hat{\varphi_1},a\hat{\varphi_2})$, thus $\text{Fix}(\hat\varphi)=\{(a,a\hat{\varphi_2})\mid a\in \text{Fix}(\hat{\varphi_1})\}$.

\begin{proposition}
Let $\varphi$ be a type II uniformly continuous endomorphism of $\mathbb Z^m \times F_n$, with $n>1$. Then, $\text{Sing}(\hat\varphi)=\{(a,a\hat{\varphi_2})\mid a\in \text{Sing}(\hat{\varphi_1})\}$. Consequentely, $\text{Reg}(\hat\varphi)=\{(a,a\hat{\varphi_2})\mid a\in \text{Reg}(\hat{\varphi_1})\}$.
\end{proposition}
\noindent \textit{Proof.} 
We start by showing that  $\text{Sing}(\hat\varphi)\subseteq \{(a,a\hat{\varphi_2})\mid a\in \text{Sing}(\hat{\varphi_1})\}$. Take some $(a,a\hat\varphi_2)\in (\text{Fix}(\varphi))^c$ with $a\in \text{Fix}(\hat\varphi_1)$. Then, for every $\varepsilon>0$, the open ball of radius $\varepsilon$ centered in $(a,a\hat\varphi_2)$ contains an element $(b_\varepsilon,b_\varepsilon\varphi_2)\in \text{Fix}(\varphi)$, with $b_\varepsilon \in \Fix(\varphi_1).$ Notice that $d(a,b_\varepsilon)\leq d((a,a\hat\varphi_2),(b_\varepsilon,b_\varepsilon\varphi_2))<\varepsilon$, thus $a\in (\Fix(\varphi_1))^c.$

For the reverse inclusion, take some $a\in \Sing(\hat\varphi_1)$. As above, we know that for every $\varepsilon>0$, there is some $b_\varepsilon\in B(a;\varepsilon)\cap \Fix(\varphi_1)$. Notice that, since $\hat\varphi_2$ is uniformly continuous, for every $\varepsilon >0$, there is some $\delta_\varepsilon$ such that, for all $a,b\in  \widehat{\mathbb Z^m}$ such that $d(a,b)<\delta_\varepsilon$, we have that $d(a\hat\varphi_2,b\hat\varphi_2)<\varepsilon$.
We want to prove that $(a,a\hat\varphi_2)\in (\Fix(\varphi))^c$, by showing that, for every $\varepsilon >0$, the ball centered in $(a,a\hat\varphi_2)$ contains a fixed point of $\varphi$. So, let $\varepsilon >0$ and consider
$\delta=\min\{\delta_{\varepsilon},\varepsilon\}$.
We have that $(b_\delta,b_\delta\hat\varphi_2)\in B((a,a\hat\varphi_2);\varepsilon)$ since, by definition of $b_\delta$, we have that $d(a,b_\delta)<\delta\leq \varepsilon$ and also, $d(a,b_\delta)<\delta_{\varepsilon}$ means that $d(a\hat\varphi_2,b_\delta \hat\varphi_2)<\varepsilon$.
\qed\\

\begin{corollary}
Let $\varphi\in \End(\mathbb Z^m\times F_n)$ be a uniformly continuous type II endomorphism. Then $\Fix (\hat\varphi)$ has $\sum\limits_{i=0}^m 2^i {m\choose i}$ $\Fix(\varphi)$-orbits.
\end{corollary}
\noindent \textit{Proof.}
For $(a,a\hat\varphi_2),(b,b\hat\varphi_2)\in \Fix(\hat\varphi)$, we have that $(a,a\hat\varphi_2)$ belongs to $(\Fix\varphi)(b,b\hat\varphi_2)$ if and only if there is some $(c,c\varphi_2) \in \Fix\varphi$ such that $(a,a\hat\varphi_2)=(c,c\varphi_2)(b,b\hat\varphi_2)$, i.e., $a$ and $b$ belong to the same orbit of $\Fix(\varphi_1)$.
\qed\\


\subsection{Finiteness Conditions on Infinite Periodic Points}
\label{subseccontinfper}
We proceed in a similar way to the case of fixed points regarding periodic points. 

We know that $\Per(\varphi)$ acts on $\Per(\hat\varphi)$ on the left, since $a\varphi^p=a$ and $b\hat\varphi^q=b$ implies that $(a+b)\hat\varphi^{pq}=a\hat\varphi^{pq}+b\hat\varphi^{pq}=a+b$ and we want to count the orbits of such action.

When we consider a type I endomorphism, we have the same issue we had in the fixed points case, in the sense that we have $\Per \hat\varphi=\Per\hat\varphi_1\times\Per\hat\phi$, which might also be uncountable.

To obtain a result for type II endomorphisms, we start as above, by dealing with the free-abelian part first. Let $\varphi\in End(\mathbb Z^m)$ defined by $a\mapsto aQ$ be a uniformly continuous endomorphism. As above, given a column $Q_j$, if $Q_j\neq 0$, we call $\lambda_j$ to its nonzero entry of column and $\alpha_j$ to its row, and if $Q_j=0$ we put $\lambda_j=0$ and $\alpha_j=1$.
Also, we will define the mapping $\psi:[m]\to [m]$ mapping $i$ to $\alpha_i$. Take $\hat\varphi:\widehat{\mathbb Z^m}\to \widehat{\mathbb Z^m} $ to be its continuous extension to the completion.

This way, we have that  $aQ=[\lambda_ia_{\alpha_i}]_{i\in [m]}$ and 
$$aQ^r=\left[ \left(\prod\limits_{j=1}^{r}\lambda_{i\psi^{j-1}}\right)a_{i\psi^{r}}\right]_{i\in[m]}.$$ 
To lighten notation, for $i\in[m]$ and $r\in \mathbb N$ we will write 
\begin{align}
\pi_i^{(r)}:=\prod\limits_{j=1}^{r}\lambda_{i\psi^{j-1}}
\end{align}

This notation will be used throughout this paper.

\begin{proposition}
\label{finitefaper}
Let $\varphi\in \End(\mathbb Z^m)$ be a uniformly continuous endomorphism. Then $\Per (\hat\varphi)$ has $\sum\limits_{i=0}^m 2^i {m\choose i}$ $\Per(\varphi)$-orbits.
\end{proposition}
\noindent \textit{Proof.}
Let $a\in \Per (\hat \varphi)$. As done in the fixed point case, we will prove that, for $b\in \Per(\hat\varphi)$, we have that $b\in (\Per\varphi)a$ if and only if their infinite entries coincide. and that suffices.

Clearly, if $b\in \Per(\hat\varphi)$ is such that $b=c+a$ for some $c\in \Per(\varphi)$, then the infinite entries of $a$ and $b$ coincide.

Now, suppose $a$ and $b$ are two infinite periodic points whose infinite entries coincide. 
Then, there are $p,q\in\mathbb N$ such that $a\hat\varphi^p=a$ and $b\hat\varphi^q=b$, so $a\hat\varphi^{pq}=a$ and $b\hat\varphi^{pq}=b$. Consider $c\in \mathbb Z^m$ defined by $c_i=0$ if $a_i,b_i\in\{+ \infty,-\infty\}$ and $c_i=b_i-a_i$ otherwise. Clearly, $b=c+a$. We only have to check that $c\in \Per(\varphi)$ and for that, we will show that $c\varphi^{pq}=c.$ We have that, for $r>0$, the mapping $\varphi^r$ is defined by 
$$[c_i]_{i\in [m]}\mapsto \left[\pi_1^{(r)}c_{i\psi^{r}}\right]_{i\in [m]}.$$
Now, we only have to see that, for every $i\in[m]$, we have that $c_i=\pi_i^{(pq)}c_{i\psi^{pq}}$. Let $i\in [m]$ such that $a_i,b_i\in\{+\infty,-\infty\}$ and so $c_i=0$. Since $|a_i|=\infty$, $a_i=(a\hat\varphi^{pq})_i=\pi_i^{(pq)}a_{i\psi^{pq}}$ and all $\lambda_k$'s are finite, then we have that $a_{i\psi^{pq}} \in\{+\infty,-\infty\}$ (and the same holds for $b_{i\psi^{pq}}$), thus, $c_{i\psi^{pq}}=0$. 
Then, we have that $\pi_i^{(pq)}c_{i\psi^{pq}}=0=c_i$. Now, take $i\in[m]$ such that $a_i,b_i\not\in \{+\infty,-\infty\}$. Then $c_i=b_i-a_i=(b\hat\varphi^{pq})_i-(a\hat\varphi^{pq})_i=\pi_i^{(pq)}b_{i\psi^{pq}}- \pi_i^{(pq)}a_{i\psi^{pq}}=\pi_i^{(pq)}(b_{i\psi^{pq}}-a_{i\psi^{pq}})=\pi_i^{(pq)}c_{i\psi^{pq}}$, since $a_{i\psi^{pq}}$ and $b_{i\psi^{pq}}$ are both finite.
\qed\\

Now, let $\varphi$ be a type II uniformly continuous endomorphism of $\mathbb Z^m \times F_n$, with $n>1$. Then $\varphi:\mathbb Z^m \times F_n \to \mathbb Z^m\times F_n$ is given by $(a,u)\mapsto (aQ,z^{\lambda a_k})$, for some $ Q\in \mathcal M_m(\mathbb Z)$ such that every column of $ Q$ contains at most one nonzero entry and $0\neq\lambda\in \mathbb Z$, $k\in [m]$. Consider $\varphi_1\in End(\mathbb Z^m)$  to be defined as $a\mapsto aQ$ and $\varphi_2: \mathbb Z^m\to F_n$ that maps $a$ to $z^{\lambda a_k}$, which are both uniformly continuous.

By uniqueness of extension, we have that $\hat\varphi: \widehat{\mathbb Z^m}\times \widehat{F_n}\to \widehat{\mathbb Z^m}\times \widehat{F_n}$ is defined by $(a,u)\mapsto (a\hat{\varphi_1},a\hat{\varphi_2})$.
Hence, if $(a,u)\in \Per(\hat\varphi)$, then $a\in \Per(\hat\varphi_1).$\\

\begin{proposition}
Let $\varphi\in \End(\mathbb Z^m\times F_n)$ be a uniformly continuous type II endomorphism. Then $\Per (\hat\varphi)$ has $\sum\limits_{i=0}^m 2^i {m\choose i}$ $\Per(\varphi)$-orbits.
\end{proposition}
\noindent \textit{Proof.}
 It is easy to see, by induction on $r$ that, for every $r>0$, we have that $(a,u)\hat\varphi^r=(a\hat\varphi_1^r,a\hat\varphi_1^{r-1}\hat\varphi_2)$. Indeed, it is true for $r=1$ and if we have that $(a,u)\hat\varphi^r=(a\hat\varphi_1^r,a\hat\varphi_1^{r-1}\hat\varphi_2)$, then $(a,u)\hat\varphi^{r+1}=(a\hat\varphi_1^r,a\hat\varphi_1^{r-1}\hat\varphi_2)\hat\varphi=(a\hat\varphi_1^{r+1},a\hat\varphi_1^{r}\hat\varphi_2).$ So we have that 
\begin{align}
\label{imgrtype2}
(a,u)\hat\varphi^r=\left([\pi_i^ra_{i\psi^r}]_{i\in[m]},z^{\lambda\pi_k^{r-1}a_{k\psi^{r-1}}}\right).
\end{align} 
Let $(a,u),(b,v)\in\Per(\hat\varphi)$. We have that $(b,v)\in (\Per \varphi)(a,u)$ if and only if $b\in (\Per(\varphi_1))a$. Indeed, if $b\in (Per(\varphi_1))a$, let $c\in \Per(\varphi_1)$ defined as in the proof of Proposition \ref{finitefaper}, such that $b=c+a$ and denote by $p,q$ the periods of $a$ and $b$, respectively. 

Then $(c,c\hat\varphi_1^{pq-1}\hat\varphi_2)\varphi^{pq}=(c,c\hat\varphi_1^{pq-1}\hat\varphi_2)$, so $(c,c\hat\varphi_1^{pq-1}\hat\varphi_2)\in \Per(\varphi)$ and 
\begin{align*}
(b,v)&=(b\hat\varphi_1^{pq},b\hat\varphi_1^{pq-1}\hat\varphi_2)=((c+a)\hat\varphi_1^{pq},(c+a)\hat\varphi_1^{pq-1}\hat\varphi_2)\\
&=(c\hat\varphi_1^{pq}+a\hat\varphi_1^{pq},c\hat\varphi_1^{pq-1}\hat\varphi_2a\hat\varphi_1^{pq-1}\hat\varphi_2)\\&=(c\hat\varphi_1^{pq},c\hat\varphi_1^{pq-1}\hat\varphi_2)(a\hat\varphi_1^{pq},a\hat\varphi_1^{pq-1}\hat\varphi_2)\\
&=(c,c\hat\varphi_1^{pq-1}\hat\varphi_2)(a,u).
\end{align*}
Hence, we have that $\Per (\hat\varphi)$ has $\sum\limits_{i=0}^m 2^i {m\choose i}$ $\Per(\varphi)$-orbits.
\qed\\

\subsection{Classification of the Infinite Fixed Points for Automorphisms}
\label{classinffixed}
Recall that a uniformly continuous automorphism of $\mathbb Z^m\times F_n$ is defined as $(a,u)\mapsto (aQ,u\phi)$, where $Q$ is a uniform matrix and $\phi\in Aut(F_n)$. 
As above, we define $\varphi_1:\mathbb Z^m\to \mathbb Z^m$ that maps $a$ to $aQ$ and we have that $\Fix(\hat\varphi)=\Fix(\hat\varphi_1)\times \Fix(\hat\phi)$

We are interested in classifying infinite fixed points. We start by presenting some standard definitions.
\begin{definition}
An infinite fixed point $\alpha\in \Fix(\hat\varphi)$ is 
\begin{itemize}
\item an \emph{attractor} if $$\exists \varepsilon >0\,\,\forall\beta\in\widehat{\mathbb Z^m \times F_n} (d(\alpha,\beta)<\varepsilon\Rightarrow \lim\limits_{n\to +\infty} \beta\hat\varphi^n=\alpha)$$
\item a \emph{repeller} if   $$\exists \varepsilon >0\,\,\forall\beta\in\widehat{\mathbb Z^m \times F_n}(d(\alpha,\beta)<\varepsilon\Rightarrow \lim\limits_{n\to +\infty} \beta\hat\varphi^{-n}=\alpha)$$
\end{itemize}
\end{definition}

We will only consider these concepts regarding automorphisms because this definition of a repeller assumes the existence of an inverse. 
It is well known that, for an automorphism of a (virtually) free group, singular fixed points cannot be attractors nor repellers  (see Proposition 1.1 in \cite{[GJLL98]}) and that every regular fixed point must be either an attractor or a repeller (see \cite{[Sil13]}). We will start by investigating what happens in the free-abelian part.

\begin{proposition}
Let $\varphi\in End(\mathbb Z^m)$ be an endomorphism defined by $a\mapsto aQ$, where $Q$ is a uniform matrix. Then $\Sing (\hat\varphi)=\Fix(\hat\varphi).$
\end{proposition}
\noindent \textit{Proof.} 
By definition, $\Sing(\hat\varphi)\subseteq \Fix(\hat\varphi).$

Let $\pi\in S_m$ be a permutation such that $\varphi$ maps $[a_i]_{i\in [m]}$ to $[\lambda_i a_{\pi(i)}]_{i\in[m]}$, and $\lambda_i=\pm 1$. Then $\Fix(\varphi)=\{a\in \mathbb Z^m \mid \forall i\in [m],\, a_i=\lambda_i a_{\pi(i)}\}$ and $\Fix(\hat\varphi)=\{a\in \widehat{\mathbb Z^m} \mid \forall i\in [m],\, a_i=\lambda_i a_{\pi(i)}\}$. Given $a\in \Fix(\hat\varphi)$ and $\varepsilon >0$, choosing some 
$$n>\max\limits_{a_i\in\mathbb Z}\left\{\left|a_i\right|,\left\lceil log_2\left(\frac 1 \varepsilon\right)\right\rceil\right\},$$
 consider $b\in \mathbb Z^m$ such that $b_i=n$ if $a_i=+\infty$; $b_i=-n$ if $a_i=-\infty$ and $b_i=a_i$, otherwise. Then $b_i\in \Fix(\varphi)$ and $d(a,b)<\varepsilon$. Thus, $b$ is a point of closure of $\Fix(\varphi)$ and we are done.
\qed\\

\begin{proposition}
Let $\varphi\in End(\mathbb Z^m)$ be defined by $a\mapsto aQ$, where $Q$ is a uniform matrix. Then an infinite fixed point $\alpha\in \Fix(\hat\varphi)\setminus \Fix(\varphi)$ is neither an attractor nor a repeller.
\end{proposition}

\noindent \textit{Proof.}
Let $a\in \Fix\hat\varphi\setminus \Fix\varphi$. Let $\varepsilon >0$ and define $q=\lceil log_2(\frac 1 \varepsilon)\rceil$. Take $p=\max\limits_{a_i\in \mathbb Z}\{q,|a_i|\}$. Take $b\in \mathbb Z^m$ such that $b_i=p$ for every $i$ such that $a_i=+\infty$; $b_i=-p$ for every $i$ such that $a_i=-\infty$ and $b_i=a_i$ otherwise. Then $d(a,b)<\varepsilon$ but $b\hat\varphi^n\not\to a$ since $\max\limits_{i\in [m]} |b_i|=p$ and applying $\hat\varphi$ simply changes order and signal of the entries, so for every $n\in \mathbb N$, we have that 
$\max\limits_{i\in [m]} \{|({b\hat\varphi^n})_i|\}=n$. Hence $d(a,b\hat\varphi^n)\geq 2^{-p}$ since there is some $k$  such that $a_k\in\{+\infty,-\infty\}$ and $d(a_k,b_k)\geq d(a_k, sgn(a_k)|p|)=2^{-p}$.
The repeller case is analogous, since the inverse of a uniform matrix is uniform.
\qed\\

So, when an endomorphism is given by a uniform matrix, no infinite fixed point is an attractor or a repeller. The next result shows how that impacts the case of a general uniformly continuous automorphism, providing a full classification of infinite fixed points.

\begin{theorem}
An infinite fixed point $\alpha=(a,u)$, where $a\in \Fix(\hat\varphi_1)$ and $u\in\Fix(\hat\phi)$ is an attractor (resp. repeller) if and only if $a$ and $u$ are attractors (resp. repellers) for $\hat\varphi_1$ and $\hat\phi$, respectively. 
\end{theorem}
\noindent \textit{Proof.} Let  $\alpha=(a,u)$ be an  infinite fixed point, where $a\in \Fix(\hat\varphi_1)$ and $u\in\Fix(\hat\phi)$. Clearly if $a\in \Fix(\hat\varphi_1)$ and $u\in \Fix(\hat\phi)$ are attractors, then, $(a,u)\in\Fix(\hat\varphi)$ is an attractor. Indeed, in that case, there are $\varepsilon_1>0$ and $\varepsilon_2>0$ such that 
$$\forall b\in\widehat{\mathbb Z^m},\, \left(d(a,b)<\varepsilon_1\Rightarrow \lim\limits_{n\to +\infty} b\hat\varphi_1^n=a\right)$$ 
and 
$$\forall v\in\widehat{ F_n},\, \left(d(u,v)<\varepsilon_2\Rightarrow \lim\limits_{n\to +\infty} v\hat\phi^n=u\right).$$
Thus, taking $\varepsilon=\min\{\varepsilon_1,\varepsilon_2\}$, we have that 
\begin{align*}
\forall (b,v)\in \widehat{\mathbb Z^m\times F_n} (d((a,u),(b,v))<\varepsilon &\Rightarrow d(a,b)<\varepsilon \wedge d(u,v)<\varepsilon \\
&\Rightarrow \lim\limits_{n\to+\infty} b\hat\varphi_1^n=a \wedge \lim\limits_{n\to+\infty} v\hat\phi^n=u\\
&\Rightarrow \lim\limits_{n\to +\infty}(b,v)\hat\varphi^n=(a,u).
\end{align*}

Conversely, suppose w.l.o.g that $a$ is not an attractor for $\hat\varphi_1$. Then, for every $\varepsilon >0$, there is some $b_\varepsilon\in \widehat{\mathbb Z^m}$ such that $d(a,b)<\varepsilon$ but $b_\varepsilon\hat\varphi_1\not\to a$. In this case, we have that, for every $\varepsilon >0$, $d((a,u),(b_\varepsilon,u))<\varepsilon$ and $(b_\varepsilon, u)\hat\varphi^n=(b_\varepsilon\hat\varphi_1^n,u)\not\to (a,u)$.
\qed\\

\begin{corollary}
\label{classinffix}
Let $\varphi\in Aut(\mathbb Z^m\times F_n)$ be a uniformly continuous automorphism such that $(a,u)\hat\varphi=(a\hat\varphi_1,b\hat\phi)$, where $\varphi_1$ is given by a uniform matrix and $\phi\in Aut(F_n)$. Then an infinite fixed point $(a,u)\in \Fix(\hat\varphi)\setminus \Fix(\varphi)$ is an attractor (resp. repeller) if and only if $a\in Fix(\varphi_1)$ and $u$ is an attractor (resp. repeller) for $\hat\phi$.
\end{corollary}

Notice that, given an infinite attractor (resp. repeller) $u\in \hat F_n$ and denoting by $S_u$ the set of points attracted (resp. repelled) to it, then the set of points attracted (resp. repelled) by $(a,u)$ is given by $T_{(a,u)}=\{(a,y)\mid y\in S_u\}$.

\section{Dynamics of infinite points}
This section is devoted to the study of the dynamics of the extension of a uniformly continuous endomorphism to the completion. We will prove that if $\varphi$ is a uniformly continuous automorphism, then the dynamics of $\hat\varphi$ is asymptotically periodic and a more general result for type II uniformly continuous endomorphisms: that every point in the completion is either periodic or wandering.
\label{dyninfpoints}
\subsection{The automorphism case}
We have that a uniformly continuous automorphism of $\mathbb Z^m\times F_n$ is defined as $(a,u)\mapsto (aQ,u\phi)$, where $Q$ is a uniform matrix and $\phi\in Aut(F_n)$. 

We start by observing that, for the abelian part, the dynamics is simple in the sense that every point is periodic.
\begin{proposition}
\label{tudoperiodico}
Let $\varphi\in Aut(\mathbb Z^m)$ be defined by $a\mapsto aQ$, where $Q$ is a uniform matrix and consider $\hat\varphi$ to be its continuous extension to the completion. Then,
there is some constant $p\leq 2^m m!$ such that $\hat\varphi^p=Id$. Hence, $\Per(\hat\varphi)=\widehat{\mathbb Z^m}$ and the period of every element divides $p$. 
\end{proposition}
\noindent \textit{Proof.}
There are only $2^mm!$ distinct uniform $m\times m$ matrices, so there are $0<p<q\leq2^mm!+1$ such that $Q^p=Q^q$, thus $I_m=Q^{q-p}$. 
\qed\\

We now present some standard dynamical definitions, that will be useful to the classification of infinite points.
\begin{definition}
Let $G$ be a group and $\varphi\in End(G).$ A point $x\in G$ is said to be a \emph{$\varphi$-wandering point} if there is a neighbourhood $U$ of $x$ and a positive integer $N$ such that for all $n>N$, we have that $U\varphi^n\cap U=\emptyset.$ When it is clear, we simply say $x$ is a wandering point.
\end{definition}

\begin{definition}
Let $G$ be a group and $\varphi\in End(G).$ A point $x\in G$ is said to be a \emph{$\varphi$-recurrent point} if, for every neighbourhood $U$ of $x$, there exists $n>0$ such that $x\varphi^n\in U$. When it is clear, we simply say $x$ is a recurrent point.
\end{definition}

\begin{definition}
Let $f$ be a homeomorphism of a compact space $K$. Given $y\in K$, the \emph{$\omega$-limit set} $\omega(y,f)$, or simply $\omega(y)$, is the set of limit points of the sequence $f^n(y)$ as $n\to + \infty$.
\end{definition}

\begin{definition}
\label{assimpperdinamic}
Let $G$ be a group. A uniformly continuous endomorphism $\varphi\in End(G)$ has \emph{asymptotically periodic dynamics} on $\hat G$ if there exists $q\geq 1$ such that, for every $x\in \hat G$, the sequence $x\hat\varphi^{qn}$ converges to a fixed point of $\hat\varphi^q$.
\end{definition}

The next proposition shows how being $\hat\varphi$-recurrent (resp. wandering) relates with being $\hat\phi$-recurrent (resp. wandering).

\begin{proposition}
\label{equivs}
Let $\varphi\in End(\mathbb Z^m \times F_n)$ be a uniformly continuous endomorphism  defined by $(a,u)\mapsto (a\varphi_1,u\phi)$, where $\varphi_1\in End(\mathbb Z^m)$  and $\phi\in End(F_n)$ and $(a,u)\in\widehat{\mathbb Z^m\times F_n}$. We have the following:
\begin{enumerate}
\item $(a,u)$ is $\hat\varphi$-periodic $\Rightarrow u$ is $\hat\phi$-periodic.
\item  $(a,u)$ is  $\hat\varphi$-wandering $\Leftarrow u$ is $\hat\phi$-wandering.
\item  $(a,u)$ is  $\hat\varphi$-recurrent $\Rightarrow u$ is $\hat\phi$-recurrent.
\end{enumerate}
\end{proposition}
\noindent \textit{Proof.} 
\begin{enumerate}
\item This is obvious. Observe that if $Q$ is uniform, which is the case when we deal with automorphisms then the reverse implication holds as well, by Propostion \ref{tudoperiodico}.
\item Suppose $u$ is $\hat\phi$-wandering. Then, take $\varepsilon>0$ and $N\in \mathbb N$ such that for all $n>N$, we have that $B(u;\varepsilon)\hat\phi^n\cap B(u;\varepsilon)=\emptyset$. Let $n>N$ and take $V=B((a,u);\varepsilon)$ and $(b,v)\in V$. Then $v\in B(u;\varepsilon)$, thus $v\hat\phi^n\not\in B(u;\varepsilon)$ and $(b,v)\hat\varphi^n=(b\hat\varphi_1^n,v\hat\phi^n)\not\in V$. Since $(b,v)$ is an arbitrary point of $V$, we have that $(a,u)$ is $\hat\phi$-wandering.
\item Suppose $(a,u)$ is $\hat\varphi$-recurrent. Take $\varepsilon>0$. There is some $n\in \mathbb N$ such that $(a,u)\hat\varphi^n\in B((a,u);\varepsilon)$ and so $u\hat\phi^n\in B(u;\varepsilon).$\qed
\end{enumerate}

Notice that,  in case $u\in F_n$ and $\phi\in End(F_n)$, if $u$ is nonwandering, then it must be periodic, since we can take taking $U=\{u\}$, we have that there is some $n$ such that $U\hat\phi^n\cap U\neq 0$ and so, it is periodic. So, in case $\phi\in Aut(\mathbb Z^m\times F_n)$, we have that a point $(a,u)\in \hat{\mathbb Z^m}\times F_n$ must be periodic or wandering, by Proposition \ref{equivs}. 

We now present two results from \cite{[LL08]} regarding free groups automorphisms, that will be very useful in this case.
 
\begin{lemma}\cite{[LL08]} Let $f$ be a homeomorphism of a compact space $K$. Given $y\in K$ and $q\geq 1$, the following conditions are equivalent:
\begin{enumerate}
\item  $\omega(y)$ is finite and has $q$ elements.
\item $\omega(y)$ is a periodic orbit of order $q$.
\item The sequence $f^{qn}(y)$ converges as $n\to +\infty$, and $q$ is minimal for this property.
\end{enumerate}
Given $p\geq 2$, the set $\omega(y,f^p)$ is finite if and only if $\omega(y,f)$ is finite.\\
\end{lemma}

If these equivalent conditions hold, we say that the point $y$ is \emph{asymptotically periodic}. If every point is asymptotically periodic, then the endomorphism has asymptotically periodic dynamics (this definition is equivalent to Definition \ref{assimpperdinamic}).

\begin{theorem}\cite{[LL08]}
Every automorphism $\alpha\in Aut(F_n)$ has asymptotically periodic dynamics. 
\end{theorem}

Using these results, we obtain the same result for free-abelian times free groups.
\begin{theorem}
\label{asympdyn}
Every uniformly continuous automorphism $\varphi\in Aut(\mathbb Z^m \times F_n)$ has asymptotically periodic dynamics on $\widehat{\mathbb Z^m\times F_n}$. 
\end{theorem}
\noindent \textit{Proof.} 
Let $\varphi$ be a uniformly continuous automorphism defined by $(a,u)\mapsto (a\varphi_1,u\phi)$, where $\varphi_1\in Aut(\mathbb Z^m)$ is given by a uniform matrix and $\phi\in Aut(F_n).$ We will prove that for every $(a,u)\in \widehat{\mathbb Z^m\times F_n}$, we have that $\omega((a,u),\varphi)$ is a periodic orbit. Take $(b,v)\in \omega((a,u),\varphi)$. Then, there is an increasing sequence $\{n_k\}$ such that $(a,u)\hat\varphi^{n_k}\to (b,v)$. This means that $a\hat\varphi_1^{n_k}\to b$ and so $b\in\omega(a,\varphi_1)$ and $u\hat\phi^{n_k}\to v$ and so $v\in  \omega(u,\phi)$ and $(b,v)$ must be a periodic point belonging to the orbit of $(a,x)$, for $x\in \omega(u,\phi)$, since $\omega(a,\hat\varphi_1)$ is the orbit of $a$. Obviously, the entire orbit must belong to $\omega((a,u),\varphi)$ since $(a,u)\hat\varphi^{n_k+p}\to(b,v)\hat\varphi^p$, for every $p\geq 0.$
 \qed

\begin{definition}
Let $\varphi$ be an endomorphism of $F_n$. We say that $\varphi$ is \emph{eventually length-nondecreasing} if $\exists p\in \mathbb N\,\forall u\in F_n\, (|u|\geq p\Rightarrow |u\varphi|\geq|u|).$
\end{definition}

We can prove that the dichotomy periodic vs wandering always holds when we deal with an eventually length-nondecreasing automorphism.
\begin{proposition}
\label{ogrande}
Let $\varphi$ be an eventually length-nondecreasing endomorphism of $F_n$. Then a point $u\in \hat{F_n}$ is either recurrent or wandering. If $\varphi$ is an automorphism, then every point is either periodic or wandering.
\end{proposition}
\noindent \textit{Proof.} 
Let $\phi$ be an eventually length-nondecreasing endomorphism of $F_n$. In particular, it is uniformly continuous (see lemma 6.2 in \cite{[CS09a]}). Take $p\in \mathbb N$ such that 
\begin{align}
\label{lengthnondec}
\forall u\in F_n\,(|u|\geq p\Rightarrow |u\varphi|\geq|u|)
\end{align}
 and $u\in \hat F_n$. Suppose $u$ is not wandering and take $\varepsilon<\frac 1 {2^p}$. Then, for every $N\in \mathbb N$ there is some $r>N$ such that $B(u;\varepsilon)\hat\phi^r\cap B(u;\varepsilon)\neq \emptyset$ and so there is some $v_r\in\hat F_n$ such that $d(u,v_r)<\varepsilon$ and $d(u,v_r\hat\phi^r)<\varepsilon$.  So, $|u\wedge v|>p$  and $|(u\wedge v)\hat\phi^r|>p$, by (\ref{lengthnondec}). Since $|u\wedge v\hat\phi^r|>p$, we have that $|u\wedge u\phi^r|>p$, hence $u\hat\phi^r\in B(u;\varepsilon)$ and $u$ is recurrent.

In case $\varphi$ is an automorphism, we have that every recurrent point must be periodic since such a point must belong to its limit set, which is a periodic orbit.
 \qed\\

Notice that this implies that if $\varphi=(\varphi_1,\phi)\in Aut(\mathbb Z^m\times F_n)$ is such that $\phi$ is length-nondecreasing, then every point $(a,u)$ is either wandering or periodic. 
\\

\subsection{Type II Endomorphisms}
The purpose of this subsection is to establish the dichotomy periodic vs wandering for Type II endomorphisms.
Recall the notation introduced in \ref{subseccontinfper}. Also, recall (\ref{imgrtype2}) and the decomposition $z=w\tilde z w^{-1}$where $\tilde z$ is the cyclically reduced core of $z$ and the definition of $B_i$ for $i\in[m]$.

\begin{remark}
Assume $n>1$. The set of periodic points of the extension of a uniformly continuous type II endomorphism to the completion is not dense in the entire space, even when we restrict ourselves to the boundary. Indeed, if we take a point $(a,u)$ such that $u$ does not share a prefix with $z$ and  $z^{-1}$, then $B\left((a,u);\frac 1 2\right)$  does not contain a periodic point. Also, the system does not admit the existence of a dense orbit: Indeed, given $(a,u)\in\widehat{\mathbb Z^m\times F_n}$, choosing a point $(b,v)\in\widehat{\mathbb Z^m\times F_n}$ such that $b\neq a$ and $v$ doesn't share a prefix of size with neither $z$ nor $z^{-1}$, we have that $B((b,v);\frac 1 2)$ does not contain any point in the orbit. 

Also, in the automorphism case, we have that the first component is always periodic, so there is not a dense orbit, even when restricted to the boundary.  
\end{remark}

We will now prove two technical lemmas that will be very useful for proving the main result.
\begin{lemma}
\label{lemaux1}
Consider a uniformly continuous endomorphism $\varphi$ of a free-abelian group $\mathbb Z^m$ and take $i\in[m]$ and some positive integer $r>m$. Then, the following conditions are equivalent:
\begin{enumerate}[(i)]
\item $\exists N \in \mathbb N\,\forall p>N \;|\lambda_{i\psi^p}|=1$
\item $\exists N \leq m \,\forall p>N \;|\lambda_{i\psi^p}|=1$
\item  $|\pi_{i\psi^{tr}}^{(r)}|=1, \text{ for every positive integer $t$}$
\item  $|\pi_{i\psi^{tr}}^{(r)}|=1, \text{ for some positive integer $t$}$
\end{enumerate}
\end{lemma}
\noindent \textit{Proof.} It is obvious that (ii) $\Rightarrow $ (i), (ii) $\Rightarrow$ (iii) and  (iii) $\Rightarrow$ (iv). We will prove that (i) $\Rightarrow$ (ii)  and that (iv)$\Rightarrow $ (i).\\

(i)$\Rightarrow$ (ii):  Suppose that there is some $N>m$ such that for every $p>N$ we have that $|\lambda_{i\psi^p}|=1$. We have that $\psi$ maps $[m]$ to a subset of $[m]$, and so, for every $i\in[m]$, there is some $k_i\leq m$ such that $i\psi^{m+1}=i\psi^{k_i}.$ This way, we have a periodic orbit (can be fixed) of $\psi$ given by $\{i\psi^{k_i},\ldots i\psi^m\}$. So, for every $p\geq N$, we define $j_p\in\{k_i,\ldots, m\}$ to be such that  $i\psi^{j_p}=i\psi^p$. Also, if for some $j>j_N$, we had $|\lambda_{i\psi^j}|>1$, we could obtain $p$ arbitrarily large such that $|\lambda_{i\psi^p}|>1$, which is absurd. So, we have that, 
 \begin{align*}
\forall p>j_N \;|\lambda_{i\psi^p}|=1
\end{align*}
and $j_N\leq m.$\\

(iv)$\Rightarrow $ (i) If, for some $i\in [m]$, we have that  $|\pi_{i\psi^{tr}}^{(r)}|=1$, for some $t\in\mathbb N$, then for every $j\in\{tr,\ldots,(t+1)r-1\}$, we have that $|\lambda_{i\psi^j}|=1$. In this case, for some $s\geq (t+1)r$, we have that $i\psi^s=i\psi^j$ for some $j\in \{tr,\ldots,(t+1)r-1\}$ and so $|\lambda_{i\psi^s}|=1$. Thus, (i) holds for $N=r$. \qed\\

\begin{lemma}
\label{stablambda}
Consider a nonwandering point $(a,u)\in \widehat{\mathbb Z^m \times F_n}$ such that $a$ has finite entries. Let  $\delta=\max\limits_{a_i\in \mathbb Z}\{|a_i|\}$ and $U=B\left((a,z^{+\infty});\frac 1 {2^\delta}\right)$. Consider a point $(b,v)\in U$ and a positive integer $r>m$ such that $(b,v)\hat\varphi^r\in U$. Then the conditions from Lemma \ref{lemaux1} hold for every index $i\in [m]$ such that $a_i\neq \pi_i^ra_{i\psi^r}$.

Moreover, if $u=z^{+\infty}$ and $|a_{k\psi^{r-1}}|<\infty$, then the conditions from Lemma \ref{lemaux1} hold when $i=k$.
\end{lemma}
\noindent \textit{Proof.} 
Consider a nonwandering point $(a,u)\in \widehat{\mathbb Z^m \times F_n}$ such that $a$ has finite entries, let  $\delta=\max\limits_{a_i\in \mathbb Z}\{|a_i|\}$ and $U=B((a,u);\frac 1 {2^\delta})$. Consider a point $(b,v)\in U$ and a positive integer $r>m$ such that $(b,v)\hat\varphi^r\in U$. So, for $i\in[m]$, we have that:
\begin{align}
\label{aifinito}
\text{if $|a_i|<\infty$, then }a_i=b_i=\pi_i^{(r)}b_{i\psi^r}
\end{align}
and
\begin{align}
\label{aiinfinito}
\text{if $|a_i|=\infty$, then }sgn(a_i)=sgn(b_i)=sgn(\pi_i^{(r)}b_{i\psi^r}) \text{ and } |b_i|,|\pi_i^{(r)}b_{i\psi^r}|>\delta.
\end{align}

If $a\hat\varphi_1^r=a$, then $a_i= \pi_i^{(r)}a_{i\psi^r}$ for every $i\in[m]$ and the first part of the lemma trivially holds.

If not, take $q\in[m]$ such that $a_q\neq \pi_q^{(r)}a_{q\psi^r}$. We start by observing that $a_q$ must be infinite since if that is not the case, then,  by (\ref{aifinito}), we have that $a_q=b_q=\pi_q^{(r)}b_{q\psi^{r}}$. If $\pi_q^{(r)}=0$, then $a_q=\pi_q^{(r)}b_{q\psi^{r}}=0=\pi_q^{(r)}a_{q\psi^{r}}$. If not, then $b_{q\psi^{r}}\leq b_q$, which means that $b_{q\psi^{r}}\leq\delta$ and by (\ref{aiinfinito}), we have that $a_{q\psi^{r}}$ is finite and  by (\ref{aifinito}), it follows that $a_{q\psi^{r}}=b_{q\psi^{r}}$, so $a_q=\pi_q^{(r)}b_{q\psi^{r}}=\pi_q^{(r)}a_{q\psi^{r}}$. Also, $a_{q\psi^{r}}$ must be finite, since, if it is infinite, then by (\ref{aiinfinito}) we have that $sgn(a_{q\psi^{r}})=sgn(b_{q\psi^{r}})$ and $sgn(a_q)=sgn(\pi_q^{(r)}b_{q\psi^r})$, thus $sgn(a_q)=sgn(\pi_q^{(r)}a_{q\psi^r})$ and that implies that $a_q=\pi_q^{(r)}a_{q\psi^{r}}$, since $\pi_q^{(r)}\neq 0$.

So, using (\ref{aifinito}) with $i=q\psi^r$, we have that 
\begin{align}
\label{form}
a_{q\psi^r}=b_{q\psi^r}=\pi_{q\psi^r}^{(r)}b_{q\psi^r\psi^r}=\pi_{q\psi^r}^{(r)}b_{q\psi^{2r}}.
\end{align}
We have that $a_q$ is infinite, so, by (\ref{aiinfinito}), we know that $\pi_q^{(r)}b_{q\psi^{r}}>\delta$. This means in particular that $a_{q\psi^{r}}\neq 0$, because otherwise we would have $b_{q\psi^{r}}=0$, by (\ref{aifinito}).

Suppose now that for every positive integer $t$, we have that  $|\pi_{q\psi^{tr}}^{(r)}|\neq 1$. If $\pi_{q\psi^r}^{(r)}=0$, then $a_{q\psi^r}=0$, which is absurd. Then, we have that  $|\pi_{q\psi^r}^{(r)}|>1$, and $|b_{q\psi^{2r}}|<|b_{q\psi^{r}}|<\infty$ and again, using (\ref{aifinito}) with $i=q\psi^{2r}$, we get that  $a_{q\psi^{2r}}=b_{q\psi^{2r}}=\pi_{q\psi^{2r}}^{(r)}b_{q\psi^{3r}}$. Again, if $\pi_{q\psi^{2r}}^{(r)}=0,$ then $a_{q\psi^{2r}}=0$ and by (\ref{aifinito}), it follows that $b_{q\psi^{2r}}=0$. From (\ref{form}), we reach a contradiction. So, we must have $|\pi_{q\psi^r}^{(r)}|>1$ and  $|b_{q\psi^{3r}}|<|b_{q\psi^{2r}}|<\infty.$ Proceeding like this, since the value of $|b_{q\psi^{pr}}|$, for $p\in \mathbb N$ cannot decrease indefinitely, then we must have that $b_{q\psi^r}=0$ and so $a_{q\psi^r}=0$, which is absurd.

So, the conditions from Lemma \ref{lemaux1} hold  when $i=q$. 

 If we have that $u=z^{+\infty}$ and $|a_{k\psi^{r-1}}|<\infty$, then by (\ref{aifinito}), we have that $a_{k\psi^{r-1}}=b_{k\psi^{r-1}}=\pi_{k\psi^{r-1}}^{(r)}b_{k\psi^{2r-1}}$. If the conditions from Lemma \ref{lemaux1}  do not hold when $i=k$, then using the same argument as above, we obtain $a_{k\psi ^{r-1}}=b_{k\psi ^{r-1}}=0$, which is absurd since $\left|z^{+\infty}\wedge z^{\lambda\pi_k^{r-1}b_{k\psi^{r-1}}}\right|>\delta$. 
 \qed\\
 
 Notice that for every $i$ for which the conditions from Lemma \ref{lemaux1} hold, there is some constant $B_i\geq 1$ such that any product of the form $\prod\limits_{j=s}^t\lambda_{i\psi^j}$, with $t\geq s$ is bounded above by $B_i$. Also, we remark that if follows from the proof that for  $q$ such that $a_q\neq \pi_q^{(r)}a_{q\psi^r}$ we must have $\pi_q^{(r)}\neq 0$, $a_q$ is infinite and $a_{q\psi^r}$ is finite.

\begin{theorem}
\label{typeiidych}
Let $\varphi\in End(\mathbb Z^m \times F_n)$ be a type II uniformly continuous endomorphism defined by $(a,u)\mapsto (aQ,z^{\lambda a_k})$, for some $k\in[m]$, $1\neq z\in F_n$, which is not a proper power and $Q$ such that $a\mapsto aQ$ is uniformly continuous. Consider $\hat\varphi$, its continuous extension to the completion. 
Then every point $(a,u)\in \widehat{\mathbb Z^m \times F_n}$ is either wandering or periodic.
\end{theorem}
\noindent \textit{Proof.} Let  $(a,u)\in \widehat{\mathbb Z^m \times F_n}$. Clearly, if $(a,u)$ is wandering, it is not periodic.
To prove the reverse inclusion, we will consider several cases.\\

\noindent\underline{\textbf{\textit{Case 1: $u\in F_n$}}}. Start by supposing that every entry in $a$ is infinite. In this case, $(a,u)$ is never periodic, so we will prove it is wandering. Take $U=B\left((a,u);\frac 1 {2^{|u|}}\right)$,  $r\in \mathbb N$ and $(b,v)\in U$. We have that, if $(b,v)\in U$, then $d(u,v)<\frac 1 {2^{|u|}}$, which means $u=v$ and for every $i \in [m]$, we have that $a_i=b_i$, or $a_ib_i>0$ and $|a_i|,|b_i|>|u|$. 

Then, we have that 
$(b,v)\hat\varphi^r=\left(\left[\pi_i^{(r)}b_{i\psi^{r}}\right]_{i\in[m]},z^{\lambda \pi_k^{(r-1)}b_{k\psi^{r-1}}}\right)$. This means that $\pi_k^{(r)}\geq 1$, because, if $\pi_k^{(r)}=0$, then $\pi_k^{(r)}b_{k\psi^{r}}=0<|u|$, which is absurd.
Since $|b_{k\psi^{r-1}}|>|u|$, then 
$$\left|z^{\lambda \pi_k^{(r-1)}b_{k\psi^{r-1}}}\right|=2|w|+\left|\lambda\pi_k^{(r-1)}b_{k\psi^{r-1}}\right||\tilde z|>|u|$$ 
and $(b,v)\hat\varphi^r\not\in U$. So, in this case, $(a,u)$ is wandering.

Now, we deal with the case where $a$ has finite entries. Suppose $(a,u)$ is not wandering. Then, for every neighbourhood $U$ of $(a,u)$, we have that $U\varphi^r\cap U\neq \emptyset$ for arbitrarily large $r$. Set $\delta=\max\limits_{a_i\in\mathbb Z}\{|a_i|,|u|\}$ and consider $U=B\left((a,u);\frac 1 {2^\delta}\right)$. We have that, if $(b,v)\in U$, then $u=v$ and  for $i\in [m]$, if $a_i$ is finite we have (\ref{aifinito}) and if $a_i$ is infinite, then we have (\ref{aiinfinito}). Take  $r>m$, $(b,v)\in U$ such that $(b,v)\varphi^r\in U.$ 
If $a\hat\varphi_1^r=a$, then if $\pi_k^{(r-1)}=0$, we have that $$u=z^{\lambda\pi_k^{(r-1)}b_{k\psi^{r-1}}}=1=z^{\lambda\pi_k^{(r-1)}a_{k\psi^{r-1}}}=a\hat\varphi_2^r.$$ 
If  $\pi_k^{(r-1)}\geq 1$, then, since $u=z^{\lambda\pi_k^{(r-1)}b_{k\psi^{r-1}}}$, we have that $|b_{k\psi^{r-1}}|<|u|\leq \delta$ and so $a_{k\psi^{r-1}}$ must be finite by (\ref{aiinfinito}). Thus, by (\ref{aifinito}), we have that $a_{k\psi^{r-1}}=b_{k\psi^{r-1}}$ and $$u=z^{\lambda\pi_k^{(r-1)}b_{k\psi^{r-1}}}=z^{\lambda\pi_k^{(r-1)}a_{k\psi^{r-1}}}=a\hat\varphi_2^r.$$ So, we have that if $a\hat\varphi_1^r=a$, then $(a,u)$ is periodic. If not, then by Lemma \ref{stablambda}, we have that the conditions from Lemma \ref{lemaux1} hold for every $i$  such that $a_i\neq \pi_i^{(r)}a_{i\psi^r}$.
Thus, the set $$X=\{j\in [m]\mid \text{the conditions from Lemma \ref{lemaux1} hold for }i=j\}$$ is nonempty.

Now, take $\tau=\max\{B_q\mid q\in X\}$ and $U'=B\left((a,u);\frac 1 {2^{\tau\delta}}\right)$. Notice that $U'\subseteq U$ and so Lemma \ref{stablambda} can be applied.
Since $(a,u)$ is nonwandering, there is some $r'>m$ and $(b',v')\in U'$ such that $(b',v')\hat\varphi^{r'}\in U'$. We will prove that $a\hat\varphi_1^{r'}=a$. Suppose not and take $q\in [m]$ such that $a_q\neq \pi_q^{(r')}a_{q\psi^{r'}}$. So $q\in X$ and from the proof of the Lemma  \ref{stablambda} it follows that $\pi_q^{(r')}\neq 0$, $a_q$ is infinite and $a_{q\psi^r}$ is finite. But then, since $(b',v')\hat\varphi^{r'}\in U'$, we must have $\pi_q^{(r')}b'_{q\psi^{r'}}>\delta\tau$, which is absurd since $\pi_q^{(r')}\leq\tau$ and $b'_{q\psi^{r'}}=a_{q\psi^{r'}}\leq \delta$.

As done above, we can check that $a\hat\varphi_2^{r'}=u$ and so $(a,u)$ is periodic. \\

\noindent\underline{\textbf{\textit{Case 2: $ u\in \partial F_n\setminus\{z^{+\infty},z^{-\infty}\}$}}}. 
In this case $(a,u)$ is never periodic, so we will prove it is wandering. Take $\delta=\max\{|z^{-\infty}\wedge u|,|z^{+\infty}\wedge u|\}$  and consider $U=B\left((a,u);\frac 1 {2^{\delta}}\right)$. Let $(b,v)\in U$. We have that $|v\wedge u|>\delta$ and for every $i\in[m]$, $a_i=b_i$ or $a_ib_i>0$ and $|a_i|,|b_i|>\delta$. So, for every $r\in \mathbb N$, we have that $(b,v)\hat\varphi^r\not\in U$, since 
$$\begin{cases}
\left|z^{\lambda \pi_k^{(r-1)}b_{k\psi^{r-1}}}\wedge u\right|=|z^{+\infty}\wedge u| \quad&\text{ if $\lambda\pi_k^{(r-1)}b_{k\psi^{r-1}}>0$}\\
0 &\text{ if $\lambda\pi_k^{(r-1)}b_{k\psi^{r-1}}=0$ }\\
\left |z^{\lambda\pi_k^{(r-1)}b_{k\psi^{r-1}}}\wedge u\right|=|z^{-\infty}\wedge u| \quad&\text{ if $\lambda \pi_k^{(r-1)}b_{k\psi^{r-1}}<0$}.
\end{cases}$$\\
\noindent\underline{\textbf{\textit{Case 3: $u\in\{z^{+\infty},z^{-\infty}\}$}}}. Suppose $(a,u)$ is not wandering and assume w.l.o.g. that $u=z^{+\infty}$. Suppose first that every entry of $a$ is infinite and consider $U=B\left((a,u),\frac 1{2^{|w|}}\right)$. Take $r>m$, $(b,v)\in U$ such that $(b,v)\varphi^r\in U.$ Denote the first letter of $\tilde z$ by $\tilde z_1$. We have that $w\tilde z_1$ is a prefix of $v$ and for every $i\in [m]$, either $a_i=b_i$ or $a_ib_i>0$ and $|a_i|,|b_i|>|w|$. 

From $(b,v)\hat\varphi^r\in U$ we deduce that $w\tilde z_1$ is a prefix of $z^{\lambda\pi_k^{(r-1)}b_{k\psi^{r-1}}}$, so $\lambda \pi_k^{(r-1)}b_{k\psi^{r-1}}>0$ and $\lambda \pi_k^{(r-1)}a_{k\psi^{r-1}}=+\infty$, since it has the same sign and every entry of $a$ is infinite. Also, $a\hat\varphi_1^r=a$ because $a_ib_i>0$ and $a_i(b\hat\varphi_1^r)_i>0$, so $\hat\varphi_1^r$ doesn't change the signs of the entries in $b$, thus it also does not change the ones in $a$. In that case, $(a,u)\hat\varphi^r=(a,u)$ and $a$ is periodic.

To complete the proof, we take $(a,z^{+\infty})$ such that $a$ has finite entries and suppose it is not wandering nor periodic.
Take $\delta=\max\limits_{a_i\in\mathbb Z}\{|a_i|\}$ and $U=B\left((a,z^{+\infty});\frac 1 {2^\delta}\right)$. Take $r>m$ and $(b,v)\in U$ such that $(b,v)\hat\varphi^r\in U$. We now consider two subcases:\\

\noindent\underline{\textbf{\textit{Subcase 3.1: $|a_{k\psi^{r-1}}|=\infty$}}}. We have that $$\left|z^{+\infty}\wedge z^{\lambda \pi_k^{(r-1)}b_{k\psi^{r-1}}}\right|>\delta$$ and so $\lambda \pi_k^{(r-1)}b_{k\psi^{r-1}}>0$. 
Thus, using (\ref{aiinfinito}) with $i=k\psi^{r-1}$, we have that $\lambda \pi_k^{(r-1)}a_{k\psi^{r-1}}>0$ and, since $|a_{k\psi^{r-1}}|=\infty$, we have that  $z^{\lambda\pi_k^{(r-1)}a_{k\psi^{r-1}}}=z^{+\infty}$. 
Since $(a,u)$ is not periodic, then we have that $a\hat\varphi_1^r=[\pi_i^{(r)}a_{i\psi^r}]_{i\in[m]}\neq a$. Take $q$ such that $a_q\neq \pi_q^{(r)}a_{q\psi^r}$. By Lemma \ref{stablambda}, we have that the conditions from Lemma \ref{lemaux1} hold when $i=q$.
Thus, the set 
$
X=\{ j\in [m] \mid \text{ the conditions from Lemma \ref{lemaux1} hold for } $i=j$ \}
$ is nonempty.

Consider  $\delta=\max\limits_{a_i\in\mathbb Z}\{|a_i|\}$, $\tau=\max\{B_q\mid q\in X\}$ and $\delta'=2|w|+\lambda\tau\sigma|\tilde z|$ and let $U'=B((a,u);\frac 1 {2^{\delta'}})$. Since $(a,u)$ is nonwandering, there is some $r'>m$ and $(b',v')\in U'$ such that $(b',v')\hat\varphi^{r'}\in U'$. Notice that $U'\subseteq U$ and so Lemma \ref{stablambda} can be applied. 
 We will prove that $(a,u)\hat\varphi^{r'}=(a,u)$, which is absurd.

We start by showing that $z^{\lambda \pi_k^{(r'-1)}a_{k\psi^{r'-1}}}=z^{+\infty}$. Indeed, if $|a_{k\psi^{r'-1}}|<\infty,$ then 
by the last statement of Lemma \ref{stablambda}, we have that $k\in X$.
Also, we know by (\ref{aifinito}) that $b'_{k\psi^{r'-1}}=a_{k\psi^{r'-1}}$, so 
$$\left|z^{\lambda\pi_k^{(r'-1)}b'_{k\psi^{r'-1}}}\right|=2|w|+|\lambda\pi_k^{(r'-1)}b'_{k\psi^{r'-1}}||\tilde z|\leq 2|w| +\lambda\tau\delta|\tilde z|=\delta',$$
which is absurd because $(b,v)\hat\varphi^r\in U'$ implies that 
\begin{align}
\label{palavraviz}
 \left|z^{+\infty}\wedge z^{\lambda\pi_k^{(r'-1)}b'_{k\psi^{r'-1}}}\right|>\delta'.
 \end{align}
 So this can never happen and so we must have $|a_{k\psi^{r'-1}}|=\infty$. Since we have (\ref{palavraviz}), it follows that $\lambda \pi_k^{(r'-1)}b'_{k\psi^{r'-1}}>0$. 
Thus, by (\ref{aiinfinito}), we have that $\lambda \pi_k^{(r'-1)}a_{k\psi^{r'-1}}>0$ and, since $|a_{k\psi^{r'-1}}|=\infty$, we have that  $z^{\lambda\pi_k^{(r'-1)}a_{k\psi^{r'-1}}}=z^{+\infty}$. 

We only have to see that $a\varphi_1^{r'}=a$. If that is not the case, then there is some $q\in[m]$ such that $a_q\neq \pi_q^{r'}a_{q\psi^{r'}}$. By Lemma \ref{stablambda}, we have that $q\in X$, which implies that $a_q$ is infinite and that $a_{q\psi^{r'}}$ is finite. It follows from $(\ref{aiinfinito})$ that $\pi_q^{(r')}b'_{q\psi^{r'}}>\delta'$, which is absurd since $\pi_q^{(r')}\leq \tau$ and, by (\ref{aifinito}), we have that  $b'_{q\psi^{r'}}=a_{q\psi^{r'}}\leq \delta$.\\
 
 \noindent\underline{\textbf{\textit{Subcase 3.2: $|a_{k\psi^{r-1}}|<\infty$}}} By Lemma \ref{stablambda}, we know that $k\in X$. Consider $\delta''=2|w|+\lambda\delta B_k|\tilde z|$ and $U''=B\left((a,u);\frac 1 {2^{\delta''}}\right)$. As usual, take $r''>m$ and $(b'',v'')\in U''$ such that $(b'',v'')\hat\varphi^{r''}\in U''$. We have that  
 $$\left|z^{+\infty}\wedge z^{\lambda\pi_k^{(r''-1)}b''_{k\psi^{r''-1}}}\right|>\delta''.$$
 But that cannot happen since $$\left|z^{\lambda\pi_k^{(r''-1)}b''_{k\psi^{r''-1}}}\right|=2|w|+|\lambda\pi_k^{(r''-1)}b''_{k\psi^{r''-1}}||\tilde z|\leq\delta''.$$
\qed\\

Since $\widehat{\mathbb Z^m\times F_n}$ is compact, then every $\omega$-limit set is nonempty. Since such a set cannot contain wandering points, then, for every point, its limit set is a periodic orbit, which means that every point is asymptotically periodic.\\

We also remark that most of these results should be  easily extended to the product of a free group with a finitely generated abelian group. So, consider $P= \mathbb Z_{p_1} \times \cdots \mathbb Z_{p_r}$ for some $r\in \mathbb N$ and powers of (not necessarily distinct) prime numbers $p_i$, for $i\in [r]$, endowed with the product metric given by taking the discrete metric in each component. This is a complete space. Take $G= F_n\times \mathbb Z^m \times P$ and $\varphi\in End(G)$. We have that a point of the form $(1,0, x_1,\ldots, x_r)$ must be mapped to a point of the same form. Indeed, setting $(1,0, x_1,\ldots, x_r)\varphi = (w,a,y_1,\ldots l_r)$, we have that $(1,0, x_1,\ldots, x_r)\varphi=(1,0, x_1,\ldots, x_r)^{p+1}\varphi$, thus $w=w^{p+1}$ and $a=(p+1)a$, so $w=1$ and $a=0$.
It follows that $(u,a,p)\varphi=((u,a)\psi_1, (u,a,p)\psi_2)$, for some $\psi_1\in End(F_n\times \mathbb Z^m)$ and $\psi_2:  F_n\times \mathbb Z^m \times P\to P$, i.e., the $P$-component has no influence in the $F_n\times \mathbb Z^m$-component of the image. 

If we want $\varphi$ to be uniformly continuous, we can see that $\varphi$ must be given by $(u,a,p)\varphi=((u,a)\psi_1, p\psi_2)$, for some $\psi_1\in End(F_n\times \mathbb Z^m)$ and $\psi_2\in  End(P)$, i.e., the $\F_n\times \mathbb Z^m$-component has no influence in the $P$-component of the image. Indeed, for every element $x$ in the basis of $F_n$, setting $(x,0,0,\ldots, 0)\varphi=(w,a, n_1,\ldots n_r)$, we have that every $n_i$ must be equal to $0$ because if that was not the case taking $\varepsilon <1$, for every delta, choosing $q$ such that $qp_1\cdots p_r>\log_{2}(\frac 1 \delta)$  we would have that $$d((x,0,\ldots 0)^{qp_1\cdots p_r},(x,0,\ldots 0)^{qp_1\cdots p_r+1})<\delta$$ but $$d((x,0,\ldots 0)^{qp_1\cdots p_r}\varphi,(x,0,\ldots 0)^{qp_1\cdots p_r+1}\varphi)=1$$ and the same happens when we consider an element in the basis of $\mathbb Z^m$.

Now observe that every point in $\psi_2$ must be periodic, since $P$ is finite. So, if the dichotomy  periodic vs wandering holds for $\psi_1$ (in particular, if $\psi_1$ is type II),
 taking a nonwandering point $(w,a,p)\in G$, we have that $(w,a)$ must be a nonwandering point of $\psi_1$, thus periodic, and $p$ is periodic and so $(w,a,p)$ is $\varphi$-periodic.

\section{Further work}
\label{furtherwk}
Although we were able to establish the dichotomy wandering/periodic for type II endomorphisms, we weren't able to go very far in answering this question for type I endomorphisms. It would be interesting to find conditions on $\varphi\in End(F_n)$ such that a type I endomorphism defined by $(a,u)\mapsto (aQ,u\phi)$ satisfies that property. If we prove that a nonwandering point is always recurrent, then, by the result in \cite{[LL08]}, we get the result for automorphisms.

Another interesting question could be obtaining conditions on the properties $\phi\in End(F_n)$ must satisfy so that the converse implications of the ones in Proposition \ref{equivs} also hold.

\section*{Acknowledgements}
The author is grateful to Pedro Silva for fruitful discussions of these topics, which  improved the paper.

The author was  supported by the grant SFRH/BD/145313/2019 funded by Funda\c c\~ao para a Ci\^encia e a Tecnologia (FCT).%


\end{document}